\documentclass[fleqn, 12pt]{amsart}
\usepackage{amssymb, amsfonts, amscd}
\usepackage{latexsym, rotating, enumerate, palatino, enumitem, subfigure, multirow}
\usepackage{xy} \xyoption{all}
\usepackage{pst-all}

\DeclareMathOperator{\id}{id}
\DeclareMathOperator{\locdim}{locdim}
\DeclareMathOperator{\ord}{ord}
\DeclareMathOperator{\ev}{ev}
\DeclareMathOperator{\im}{im}
\DeclareMathOperator{\core}{core}
\DeclareMathOperator{\diam}{diam}

\newcommand{\e}{\varepsilon}

\newcommand{\KK}{{\mathbb{K}}}
\newcommand{\RR}{{\mathbb{R}}}
\newcommand{\NN}{{\mathbb{N}}}
\newcommand{\CC}{{\mathbb{C}}}
\newcommand{\Toep}{{\mathcal{T}}}
\newcommand{\cG}{{\mathcal{G}}}
\newcommand{\cpctnPt}[1]{\alpha\! #1}
\newcommand{\Cs}{{$C^*$-algebra}}
\newcommand{\termDef}[1]{\textbf{#1}\index{#1}}
\newcommand{\impliesStep}[2]{''(#1) $\Rightarrow$ (#2)''}

\newtheorem{lma}{Lemma}[section]
\newtheorem{cor}[lma]{Corollary}
\newtheorem{thm}[lma]{Theorem}
\newtheorem{prop}[lma]{Proposition}
\newtheorem{defn}[lma]{Definition}
\theoremstyle{definition}
\newtheorem{pargr}[lma]{}
\newtheorem{rmk}[lma]{Remark}
\newtheorem{quest}[lma]{Question}

\begin{document}

\title[Semiprojectivity of commutative $C^*$-algebras.]{A characterization of semiprojectivity for commutative $C^*$-algebras.}%
\author{Adam P. W. S{\o}rensen and Hannes Thiel}

\address{Department of Mathematical Sciences, University of Copenhagen, Universitetsparken 5, DK-2100, Copenhagen \O, Denmark}
\email{apws@math.ku.dk}
\address{Department of Mathematical Sciences, University of Copenhagen, Universitetsparken 5, DK-2100, Copenhagen \O, Denmark}
\email{thiel@math.ku.dk}

\thanks{
This research was supported by the Danish National Research Foundation through the Centre for Symmetry and Deformation.
The second named author was partially supported by the Marie Curie Research Training Network EU-NCG
}

\subjclass[2000]%
{Primary
46L05, % General theory of C*-algebras
54C55, % Absolute neighborhood retracts and extensors
55M15, % Absolute neighborhood retracts
54F50 % Spaces of dimension $\le 1$
 ; Secondary
46L80, % K-theory and operator algebras
%46L35, % Classifications and factors of C*-algebras
%19K14, % $K\sb 0$ as an ordered group
%46L85, % Noncommutative topology
46M10, % Projective / injective objects in categories of topol. linear spaces
54F15, % Continua and generalizations
54D35, % Compactifications
54C56, % Shape theory
55P55 % Shape theory
}

\date{\today}

\begin{abstract}
    Given a compact, metric space $X$, we show that the commutative \Cs{} $C(X)$ is semiprojective if and only if $X$ is an absolute neighborhood retract of dimension at most one.
    This confirms a conjecture of Blackadar.

    Generalizing to the non-unital setting, we derive a characterization of semiprojectivity for separable, commutative \Cs{s}.
    As further application of our findings we verify two conjectures of Loring and Blackadar in the commutative case, and we give a partial answer to the question, when a commutative \Cs{} is weakly (semi-)projective.
\end{abstract}

\maketitle

%################################################################################
%################################################################################
\section{Introduction}
\label{sect:01:introduction}

\noindent
    Shape theory is a machinery that allows to focus on the global properties of a space by abstracting from its local behavior.
    This is done by approximating the space by a system of nicer spaces, and then studying this approximating system instead of the original space.
    After this idea was successfully applied to commutative spaces, it was first introduced to the noncommutative world by Effros and Kaminker, \cite{EffKam1986}.
    Soon after, noncommutative shape theory was developed to its modern form by Blackadar, \cite{Bla1985}.

    In classical shape theory one approximates a space by absolute neighborhood retracts (ANRs).
    In the noncommutative world, the role of these nice spaces is played by the semiprojective \Cs{s}.
    It is however not true that every (compact) ANR $X$ gives a semiprojective \Cs{} $C(X)$.
    In fact, already the two-disc $D^2$ is a counterexample (see \ref{prop:S03:Disc_not_SP} and \ref{pargr:S03:spaces_containing_disc}).
    This hints to a possible problem in noncommutative shape theory:
    While it easy to show that there are enough ANRs to approximate every compact metric space, the analogue for \Cs{} is not obvious at all.
    In fact it is still an open problem whether every separable \Cs{} can be written as an inductive limit of semiprojective \Cs{s}.
    Some progress on this problem was recently made by Loring and Shulman, \cite{LorShu2010}.

    Hence, it is important to know which \Cs{s} are semiprojective.
    And although semiprojectivity was modeled on ANRs, the first large class of \Cs{s} shown to be semiprojective were the highly noncommutative Cuntz-Krieger algebras, see \cite{Bla1985}.
    Since then, these results have been extended to cover all UCT Kirchberg algebras with finitely generated K-theory and free $K_1$-group, see \cite{Szy2002} and \cite{Spi2009}, and it is conjectured that in fact all Kirchberg algebras with finitely generated K-theory are semiprojective.

    Yet, the following natural question remained unanswered:

    \begin{quest} \label{quest:S01:Main_quest}
        Which commutative \Cs{s} are semiprojective?
    \end{quest}

    An important partial answer was obtained by Loring, \cite[Proposition 16.2.1, p.125]{Lor1997}, who showed that all one-dimensional CW-complexes give rise to semiprojective \Cs{s}.
    In \cite{ELP1998} this was extended to the class of one-dimensional NCCW-complexes.

    In another direction, Chigogidze and Dranishnikov recently gave a characterization of the commutative \Cs{s} that are projective: They show in \cite[Theorem 4.3]{ChiDra2010} that $C(X)$ is projective in $\mathcal{S}_1$ (the category of unital, separable \Cs{s} with unital $\ast$-homomorphisms) if and only if $X$ is an AR and $\dim(X)\leq 1$.
    Inspired by their results we obtain the following answer to question \ref{quest:S01:Main_quest}:

\begin{thm} \label{MainTheorem}
    Let $X$ be a compact, metric space.
    Then the following are equivalent:
    \begin{enumerate}[label=(\Roman*)]
	   \item $C(X)$ is semiprojective.
	   \item $X$ is an ANR and $\dim(X)\leq 1$.
    \end{enumerate}
\end{thm}

    This confirms a conjecture of Blackadar, \cite[II.8.3.8, p.163]{Bla2006}.
    We proceed as follows:

\tableofcontents

\noindent
    In section $2$ (Preliminaries), we recall the basic concepts of commutative and noncommutative shape theory, in particular the notion of an ANR and of semiprojectivity.
    \\

\noindent
    In section $3$ (Necessity), we show the implication ''(I) $\Rightarrow$ (II)'' of our main result \ref{MainTheorem}.
    The idea is to use the topological properties of higher dimensional spaces, to show that if $C(X)$ was semiprojective and $X$ an ANR of dimension at least $2$ then we could solve a lifting problem known to be unsolvable.
    \\

\noindent
    In section $4$ we study the structure of compact, one-dimensional ANRs.
    We characterize when a one-dimensional Peano continuum $X$ is an ANR, see \ref{prop:S04:TFAE_1D_Peano_ANR}.
    As it turns out, one criterium is that $X$ contains a finite subgraph that contains all homotopy information, a (homotopy) core, see \ref{prop:S04:Existence_of_core}.
    This is also equivalent to $K^\ast(X)$ being finitely generated, which is a recurring property in connection with semiprojectivity.

    The main result of this section is theorem \ref{prop:S04:Structure_1D_Peano} which describes the internal structure of a compact, one-dimensional ANR $X$.
    Starting with the homotopy core $Y_1\subset X$ there is an increasing sequence of subgraphs $Y_1\subset Y_2\subset\ldots\subset X$ that exhaust $X$, and such that $Y_{k+1}$ is obtained from $Y_k$ by simply attaching a line segment at one end to a point in $Y_k$.
    This generalizes the classical structure theorem for dendrites (which are precisely the \emph{contractible}, compact, one-dimensional ANRs).
    \\

\noindent
    In section $5$ (Sufficiency) we show the implication ''(II) $\Rightarrow$ (I)'' of \ref{MainTheorem}.
    Using the structure theorem \ref{prop:S04:Structure_1D_Peano} for $X$, we obtain subgraphs $Y_k\subset X$ such that $X\cong\varprojlim Y_k$.
    The first graph $Y_1$ contains all K-theory information, and the subsequent graphs are obtained by attaching line segments.
    Dualizing, we can write $C(X)$ as an inductive limit, $C(X) = \varinjlim C(Y_k)$.
    Since the maps $Y_{k+1}\to Y_k$ are retractions, the dual bonding morphisms $C(Y_k)\to C(Y_{k+1})$ are accessible for lifting problems.

    The main result of this section is \ref{InductiveLimitProjective}.
    Given a lifting problem $C(X)\to C/\overline{\bigcup_k J_k}$ and an initial lift from $C(Y_1)$ to some $C/J_l$, there exists a lifting from any $C(Y_k)$ to the same height, and finally a lift from the inductive limit $C(X)$ to $C/J_l$.
    This idea is central in \cite{ChiDra2010}, but it has also been used before, for instance by Blackadar in order to prove that the Cuntz algebra $\mathcal{O}_\infty$ is semiprojective.
    We note that some form of inductive limit argument seems necessary for lifting an infinite number of generators.
    We also wish to point out that Chigogidze and Dranishnikov only needed semiprojectivity, and not projectivity, in many steps of their proofs.

    The proof ''(II) $\Rightarrow$ (I)'' follows from \ref{InductiveLimitProjective} if we can find an initial lift from $C(Y_1)$.
    For this we use Loring's deep result, \cite{Lor1997}, which says that $C(Y)$ is semiprojective for every finite graph $Y$.
    We also need Loring's result to write the algebras $C(Y_k)$ as universal \Cs{s}.

    To summarize, the proof proceeds in two steps.
    First, we construct an initial lift $C(Y_1)\to C/J_l$ from the homotopy core.
    This will lift all K-theory information of $X$.
    %, and we can not lift ''all the way'' to $C$, but only to a certain level $C/J_l$.
    But once the K-theory information is lifted, we do not need to ''sink to a lower level''.
    \\

\noindent
    In section $6$ we give applications of our main result \ref{MainTheorem}.
    First, we analyze the structure of non-compact, one-dimensional ANRs.
    We give a characterization when the one-point compactification of such spaces is again an ANR, see \ref{prop:S06:cpctfn_ANR}.
    This is motivated by the fact that a \Cs{} $A$ is semiprojective if and only if its minimal unitalization $\widetilde{A}$ is semiprojective.
    For commutative \Cs{s}, the minimal unitalization corresponds to taking the one-point compactification of the underlying commutative space.
    Using the characterization of semiprojectivity for unital, separable, commutative \Cs{s} given in \ref{MainTheorem}, we derive a characterization of semiprojectivity for non-unital, separable, commutative \Cs{s}, see \ref{prop:S06:SP_for_non-compact}.

    In \ref{prop:S06:cpctfn_ANR} we also note that the one-point compactification of the considered spaces is an ANR if and only every finite-point compactification is an ANR.
    This allows us to study short exact sequences
    \begin{center}
    	\makebox{
        \xymatrix{
        0\ar[r] & I \ar[r] & A \ar[r] & F \ar[r] & 0
        \\
        }}
    \end{center}
    with $F$ finite-dimensional.
    It was conjectured by Loring and also by Blackadar, \cite[Conjecture 4.5]{Bla2004}, that in this situation $A$ is semiprojective if and only if $I$ is.
    One implication was recently proven by Dominic Enders, \cite{EndPrivat}, who showed that semiprojectivity passes to ideals when the quotient is finite-dimensional.
    The converse implication is in general not even known for $F=\CC$.
    However, in \ref{prop:S06:ideal_with_fd_quotient} we verify this conjecture under the additional assumption that $A$ is commutative.

    Then, we will study the semiprojectivity of \Cs{s} of the form $C_0(X,M_k)$.
    We derive in \ref{matrixMain} that for a separable, commutative \Cs{} $A$, the algebra $A\otimes M_k$ is semiprojective if and only if $A$ is semiprojective.
    Again, this question can be asked in general.
    It is known that semiprojectivity of $A$ implies that $A\otimes M_k$ is semiprojective as well, see \cite[Corollary 2.28]{Bla1985} and \cite[Thoerem 14.2.2, p.110]{Lor1997}.
    For the converse, it is known that semiprojectivity passes to full corners, \cite[Proposition 2.27]{Bla1985}.
    It was conjectured by Blackadar, \cite[Conjecture 4.4]{Bla2004}, that the same holds for full hereditary sub-\Cs{s}.
    Note that $A$ always is a full hereditary sub-\Cs{} of $A\otimes M_k$.
    Thus, we verify the conjecture for commutative \Cs{s}.

    As a final application, we consider the following variant of question \ref{quest:S01:Main_quest}:
    When is a commutative \Cs{} weakly (semi-)projective?
    In order to study this problem, we analyze the structure of one-dimensional approximative absolute (neighborhood) retracts, abbreviated AA(N)R.
    In \ref{prop:S07:TFAE_1D_AANR} we show that such spaces are approximated from within by finite trees (finite graphs).
    Since finite trees (finite graphs) give (semi-)projective \Cs{s}, we derive in \ref{prop:S07:1D_AANR_implies_wSP} that $C(X)$ is weakly (semi-)projective in $\mathcal{S}_1$ if $X$ is a one-dimensional AA(N)R.

    Summarizing our results, \ref{MainTheorem} and \ref{prop:S07:1D_AANR_implies_wSP}, and the result of Chigogidze and Dranishnikov, \cite[Theorem 4.3]{ChiDra2010}, we get:

\begin{thm} \label{summaryThm}
    Let $X$ be a compact, metric space with $\dim(X)\leq 1$.
    Then:

    \begin{tabular}{lll}
        (1)\quad & $C(X)$ is projective in $\mathcal{S}_1$
            & $\Leftrightarrow$ $X$ is an AR \\
        (2) & $C(X)$ is weakly projective in $\mathcal{S}_1$
            & $\Leftrightarrow$ $X$ is an AAR \\
        (3) & $C(X)$ is semiprojective $\mathcal{S}_1$
            & $\Leftrightarrow$ $X$ is an ANR \\
        (4) & $C(X)$ is weakly semiprojective $\mathcal{S}_1$
            & $\Leftrightarrow$ $X$ is an AANR \\
    \end{tabular}

\noindent
    Moreover, $C(X)$ projective or semiprojective already implies $\dim(X)\leq 1$.
    \\
\end{thm}

%################################################################################
%################################################################################
\section{Preliminaries}
\label{sect:02:preliminaries}

\noindent
    By $A,B,C,D$ we mostly denote \Cs{s}, usually assumed to be separable here, and by a morphism between \Cs{s} we understand a $\ast$-homomorphism.
    By an ideal in a \Cs{} we mean a closed, two-sided ideal.
    If $A$ is a \Cs{}, then we denote by $\widetilde{A}$ its minimal unitalization, and by $A^+$ the forced unitalization.
    Thus, if $A$ is unital, then $\widetilde{A}=A$ and $A^+\cong A\oplus\CC$.
    We use the symbol $\simeq$ to denote homotopy equivalence.

    By a map between two topological spaces we mean a continuous map.
    Given $\varepsilon>0$ and subsets $F,G\subset X$ of a metric space, we say $F$ is \termDef{$\varepsilon$-contained} in $G$, denoted by $F\subset_\varepsilon G$, if for every $x\in F$ there exists some $y\in G$ such that $d_X(x,y)<\varepsilon$.
    Given two maps $\varphi,\psi\colon X\to Y$ between metric spaces and a subset $F\subset X$ we say ''$\varphi$ and $\psi$ agree on $F$'', denoted $\varphi=^F\psi$, if $\varphi(x)=\psi(x)$ for all $x\in F$.
    If moreover $\varepsilon>0$ is given, then we say ''$\varphi$ and $\psi$ agree up to $\varepsilon$'', denoted $\varphi=_\varepsilon\psi$, if $d_Y(\varphi(x),\psi(x))<\varepsilon$ for all $x\in X$ (for normed spaces, this is usually denoted by $\|\varphi-\psi\|_\infty<\varepsilon$).
    We say ''$\varphi$ and $\psi$ agree on $F$ up to $\varepsilon$'', denoted $\varphi=_\varepsilon^F\psi$, if $d_Y(\varphi(x),\psi(x))<\varepsilon$ for all $x\in F$.

%--------------------------------------------------------------------------------
\begin{pargr}[(Approximative) absolute (neighborhood) retracts]
\label{pargr:S02:AANR}
    A metric space $X$ is an \termDef{(approximative) absolute retract}, abbreviated by \termDef{(A)AR}, if for all
pairs\footnote{A $(Y,Z)$ pair of spaces is simply a space $Y$ with a \emph{closed} subspace $Z\subset Y$.}
    $(Y,Z)$ of metric spaces and maps $f\colon Z\to X$ (and $\varepsilon>0$)
    there exists a map $g\colon Z\to X$ such that $f=g\circ\iota$ (resp. $f=_\varepsilon g\circ\iota$), where $\iota\colon Z\hookrightarrow Y$ is the inclusion map.
    This means that the following diagram can be completed to commute (up to $\varepsilon$ ):
    \begin{center}
    	\makebox{
        \xymatrix@M+=5pt{
        & Y \ar@{.>}[dl]_{g} \\
        X & Z \ar[l]^{f} \ar@{^{(}->}[u]_{\iota}
        \\
        }}
    \end{center}

    A metric space $X$ is an \termDef{(approximative) absolute neighborhood retract}, abbreviated by \termDef{(A)ANR}, if for all pairs $(Y,Z)$ of metric spaces and maps $f\colon Z\to X$ (and $\varepsilon>0$)
    there exists a neighborhood $V$ of $Z$ and a map $g\colon V\to X$ such that $f=g\circ\iota$ (resp. $f=_\varepsilon g\circ\iota$) where $\iota\colon Z\hookrightarrow V$ is the inclusion map.
    This means that the following diagram can be completed to commute (up to $\varepsilon$ ):
    \begin{center}
    	\makebox{
        \xymatrix@M+=5pt{
        & Y  \\
        & V \ar@{^{(}->}[u] \ar@{..>}[dl]_{g} \\
        X & Z \ar[l]^{f} \ar@{^{(}->}[u]_{\iota}
        }}
    \end{center}

    For details about ARs and ANRs see \cite{Bor1967}.
    We will only consider compact AARs and AANRs in this paper, and the reader is referred to \cite{Cla1971} for more details.
    \\
\end{pargr}

%--------------------------------------------------------------------------------
\noindent
    We consider shape theory for separable \Cs{s} as developed by Blackadar, \cite{Bla1985}.
    Let us shortly recall the main notions and results:
    \\

%--------------------------------------------------------------------------------
\begin{pargr}[(Weakly) (semi-)projective \Cs{s}]
\label{pargr:S02:wSP}
    Let $\mathcal{D}$ be a subcategory of the category of \Cs{s}, closed under
quotients\footnote{This means the following: Assume $B$ is a quotient \Cs{} of $A$ with quotient morphism $\pi\colon A\to B$.
If $A\in\mathcal{D}$, then $B\in\mathcal{D}$ and $\pi$ is a $\mathcal{D}$-morphism.}.
    A $\mathcal{D}$-morphism $\varphi\colon A\to B$ is called \termDef{(weakly) projective in $\mathcal{D}$} if for any \Cs{} $C$ in $\mathcal{D}$ and $\mathcal{D}$-morphism $\sigma\colon B\to C/J$ to some quotient (and finite subset $F\subset A$, $\varepsilon>0$), there exists a $\mathcal{D}$-morphism $\bar{\sigma}\colon A\to C$ such that $\pi\circ\bar{\sigma}=\sigma\circ\varphi$ (resp. $\pi\circ\bar{\sigma}=_\varepsilon^F\sigma\circ\varphi$), where $\pi\colon C\to C/J$ is the quotient morphism.
    This means that the following diagram can be completed to commute (up to $\varepsilon$ on $F$):
    \begin{center}
    	\makebox{
        \xymatrix{
        & & C \ar[d]^{\pi} \\
        A \ar[r]_{\varphi} \ar@{..>}[urr]^{\bar{\sigma}} & B \ar[r]_{\sigma}& C/J
        \\
        }}
    \end{center}
    A \Cs{} $A$ is called \termDef{(weakly) projective} in $\mathcal{D}$ if the identity morphism $\id_A\colon A\to A$ is (weakly) projective.

    A $\mathcal{D}$-morphism $\varphi\colon A\to B$ is called \termDef{(weakly) semiprojective in $\mathcal{D}$} if for any \Cs{} $C$ in $\mathcal{D}$ and increasing sequence of ideals $J_1\lhd J_2\lhd\ldots\lhd C$ and $\mathcal{D}$-morphism $\sigma\colon B\to C/\overline{\bigcup_k J_k}$ (and finite subset $F\subset A$, $\varepsilon>0$), there exists an index $k$ and a $\mathcal{D}$-morphism $\bar{\sigma}\colon A\to C/J_k$ such that $\pi_k\circ\bar{\sigma}=\sigma\circ\varphi$ (resp. $\pi_k\circ\bar{\sigma}=_\varepsilon^F\sigma\circ\varphi$), where $\pi_k\colon C/J_k\to C/\overline{\bigcup_k J_k}$ is the quotient morphism.
    This means that the following diagram can be completed to commute (up to $\varepsilon$ on $F$):
    \begin{center}
    	\makebox{
        \xymatrix{
        & & C \ar[d] \\
        & & C/J_k \ar[d]^{\pi} \\
        A \ar[r]_{\varphi} \ar@{..>}[urr]^{\psi} & B \ar[r]_{\sigma}  & C/\overline{\bigcup_k J_k}
        \\
        }}
    \end{center}
    A \Cs{} $A$ is called \termDef{(weakly) semiprojective} in $\mathcal{D}$ if the identity morphism $\id_A\colon A\to A$ is (weakly) semiprojective.

    It is well known that if $A$ is separable then $A$ is semiprojective in the category of all \Cs{s} if and only if it is in the category of separable \Cs{s}.
    If $\mathcal{D}$ is the category $\mathcal{S}$ of all separable \Cs{s} (with all $\ast$-homomorphisms), then one drops the reference to $\mathcal{D}$ and simply speaks of (weakly) (semi-)projective \Cs{s}.
    Besides $\mathcal{S}$ one often considers the category $\mathcal{S}_1$ of all \emph{unital} separable \Cs{s} with \emph{unital} $\ast$-homomorphisms as morphisms.

    A projective \Cs{} cannot have a unit.
    For a (separable) \Cs{s} $A$ we get from \cite[Proposition 2.5]{Bla1985}, see also \cite[Theorem 10.1.9, p.75]{Lor1997}, that the following are equivalent:
    \begin{enumerate}[label=(\arabic*)]
        \item $A$ is projective
        \item $\widetilde{A}$ is projective in $\mathcal{S}_1$
    \end{enumerate}

    The situation for semiprojectivity is even easier.
    A unital \Cs{} is semiprojective if and only if it is semiprojective in $\mathcal{S}_1$.
    Further, for a separable \Cs{} $A$ we get from \cite[Corollary 2.16]{Bla1985}, see also \cite[Theorem 14.1.7, p.108]{Lor1997}, that the following are equivalent:
    \begin{enumerate}[label=(\arabic*)]
        \item $A$ is semiprojective
        \item $\widetilde{A}$ is semiprojective
        \item $\widetilde{A}$ is semiprojective in $\mathcal{S}_1$
        \\
    \end{enumerate}
\end{pargr}

%--------------------------------------------------------------------------------
\begin{pargr}[Connection between (approximative) absolute (neighborhood) retracts and (weakly) (semi-)projective \Cs{s}]
\label{pargr:S02:Connection}
    Let $\mathcal{SC}$ be the full subcategory of $\mathcal{S}$ consisting of (separable) commutative \Cs{s}, and similarly let $\mathcal{SC}_1$ be the full subcategory of $\mathcal{S}_1$ consisting of (separable, unital) commutative \Cs{s}.

    In general, for a \Cs{} it is easier to be (weakly) (semi-)projective in a smaller full subcategory, since there are fewer quotients to map into.
    In particular, if a commutative \Cs{} is (weakly) (semi-)projective, then it will be (weakly) (semi-)projective with respect to $\mathcal{SC}$.
    If one compares the definitions carefully, then one gets the following equivalences for a \emph{compact}, metric space $X$ (see \cite[Proposition 2.11]{Bla1985}):

    \begin{tabular}{lll}
        (1)\quad & $C(X)$ is projective in $\mathcal{SC}_1$
            & $\Leftrightarrow$ $X$ is an AR \\
        (2) & $C(X)$ is weakly projective in $\mathcal{SC}_1$
            & $\Leftrightarrow$ $X$ is an AAR \\
        (3) & $C(X)$ is semiprojective in $\mathcal{SC}_1$
            & $\Leftrightarrow$ $X$ is an ANR \\
        (4) & $C(X)$ is weakly semiprojective in $\mathcal{SC}_1$
            & $\Leftrightarrow$ $X$ is an AANR \\
    \end{tabular}

    Thus, the notion of (weak) (semi-)projectively is a translation of the concept of an (approximate) absolute (neighborhood) retract to the world of noncommutative spaces.
    Let us clearly state a point which is used in the proof of the main theorem:
    If $C(X)$ is (weakly) (semi-)projective in $\mathcal{SC}_1$, then $X$ is an (approximate) absolute (neighborhood) retract.
    As we will see, the converse is not true in general.
    We need an assumption on the dimension of $X$.
    \\
\end{pargr}

%--------------------------------------------------------------------------------
\begin{pargr}[Covering dimension]
\label{pargr:S02:Dimension}
	By $\dim(X)$ we denote the covering dimension of a space $X$.
	By definition, $\dim(X)\leq n$ if every finite open cover $\mathcal{U}$ of $X$ can be refined by a finite open cover $\mathcal{V}$ of $X$ such that $\ord(\mathcal{V})\leq n+1$.
	Here $\ord(\mathcal{V})$ is the largest number $k$ such that there exists some point $x\in X$ that is contained in $k$ different elements of $\mathcal{V}$.

    To an open cover $\mathcal{V}$ one can naturally assign an abstract simplicial
complex\footnote{An abstract simplicial complex over a set $S$ is a family $C$ of finite subsets of $S$ such that $X\subset Y\in C$ implies $X\in C$. An element $X\in C$ with $n+1$ elements is called an $n$-simplex (of the abstract simplicial complex).}
    $\mathcal{N}(\mathcal{V})$, called the nerve of the covering.
    It is is defined as the family of finite subsets $\mathcal{V}'\subset\mathcal{V}$ with non-empty intersection, in symbols:
    \begin{align*}
        \mathcal{N}(\mathcal{V}):=\{ \mathcal{V}'\subset\mathcal{V} \text{ finite }\ :\  \bigcap\mathcal{V}'\neq\emptyset \}.
    \end{align*}
    A $n$-simplex of $\mathcal{N}(\mathcal{V})$ corresponds to a choice of $n$ different elements in the cover that have non-empty intersection.
    Given an abstract simplicial complex $C$, one can naturally associate to it a space $|C|$, called the geometric realization of $C$.
    The space $|C|$ is a polyhedron, in particular it is a CW-complex.

    Note that $\ord(\mathcal{V})\leq n+1$ if and only if the nerve $\mathcal{N}(\mathcal{V})$ of the covering $\mathcal{V}$ is an abstract simplicial set of
dimension\footnote{The dimension of an abstract simplicial set is the largest integer $k$ such that it contains a $k$-simplex.}
    $\leq n$, or equivalently the geometric realization of $|\mathcal{N}(\mathcal{V})|$ is a polyhedron of covering
dimension\footnote{The covering dimension of polyhedra, or more generally CW-complexes, is easily understood. These spaces are successively build by attaching cells of higher and higher dimension.
The (covering) dimension of a CW-complex is simply the highest dimension of a cell that was attached when building the complex.}
    $\leq n$.

    Let $\mathcal{U}$ be a finite open covering of a space $X$, and $\{e_u\ :\ U\in\mathcal{U}\}$ a partition of unity that is subordinate to $\mathcal{U}$.
    This naturally defines a map $\alpha\colon X\to|\mathcal{N}(\mathcal{U})|$ sending a point $x\in X$ to the (unique) point $\alpha(x)\in|\mathcal{N}(\mathcal{U})|$ that has ''coordinates'' $e_U(x)$.
	
	By $\locdim(X)$ we denote the local covering dimension of a space $X$.
	By definition $\locdim(X)\leq n$ if every point $x\in X$ has a closed neighborhood $D$ such that $\dim(D)\leq n$.
	If $X$ is paracompact (e.g. if it is compact, or locally compact and $\sigma$-compact), then $\locdim(X)=\dim(X)$.

	See \cite{Nag1970} for more details on nerves, polyhedra and the (local) covering dimension of a space.
    \\
\end{pargr}

%--------------------------------------------------------------------------------
\noindent
    A particularly nice class of
one-dimensional\footnote{We say a space is one-dimensional if $\dim(X)\leq 1$. So, although it sounds weird, a one-dimensional space can also be zero-dimensional. It would probably be more precise to speak of ''at most one-dimensional'' space, however the usage of the term ''one-dimensional space'' is well established.} spaces are the so-called dendrites.
    Before we look at them, let us recall some notions from continuum theory.
    A good reference is Nadler's book, \cite{Nad1992}.

    A \termDef{continuum} is a compact, connected, metric space, and a \termDef{generalized continuum} is a locally compact, connected, metric space.
    A \termDef{Peano continuum} is a locally connected continuum, and a \termDef{generalized Peano continuum} is a locally connected generalized continuum.
    By a \termDef{finite graph} we mean a graph with finitely many vertices and edges, or equivalently a compact, one-dimensional CW-complex.
    By a \termDef{finite tree} we mean a contractible finite graph.
    \\

%--------------------------------------------------------------------------------
\begin{pargr}[Dendrites]
\label{pargr:S02:Dendrites}
    A \termDef{dendrite} is a Peano continuum that does not contains a simple closed curve (i.e., there is no embedding of the circle $S^1$ into it).
    There are many other characterizations of a dendrite.
    We collect a few and we will use them without further mentioning.

    Let $X$ be a Peano continuum.
    Then $X$ is a dendrite if and only if one (or equivalently all) of the following conditions holds:
    \begin{enumerate}[label=(\arabic*)]
        \item
        $X$ is one-dimensional and contractible
        \item
        $X$ is
tree-like\footnote{A (compact, metric) space $X$ is tree-like, if for every $\e>0$ there exists a finite tree $T$ and a map $f\colon X\to T$ onto $T$ such that $\diam(f^{-1}(y))<\e$ for all $y\in T$.}.
        \item
        $X$ is
dendritic\footnote{A space $X$ is called dendritic, if any two points of $X$ can be separated by the omission of a third point}
        \item
        $X$ is hereditarily
unicoherent\footnote{A continuum $X$ is called unicoherent if for each two subcontinua $Y_1,Y_2\subset X$ with $X=Y_1\cup Y_2$ the intersection $Y_1\cap Y_2$ is a continuum (i.e. connected).
A continuum is called hereditarily unicoherent if all its subcontinua are unicoherent.}.
    \end{enumerate}
    For more information about dendrites see \cite[Chapter 10]{Nad1992}, \cite{Lel1976}, \cite{CasCha1960}.
    \\
\end{pargr}

%################################################################################
%################################################################################
\section{One implication of the main theorem: Necessity}
\label{sect:03:necessity}

%--------------------------------------------------------------------------------
\begin{prop}
\label{prop:S03:SP_gives_1D}
	Let $C(X)$ be a unital, separable \Cs{} that is semiprojective.
	Then $X$ is a compact ANR with $\dim(X)\leq 1$.
\end{prop}
\begin{proof}
    Assume such a $C(X)$ is given.
    Then $X$ is a compact, metric space.
    As noted in \ref{pargr:S02:Connection}, semiprojectivity (in $\mathcal{S}_1$) implies semiprojectivity in the full subcategory $\mathcal{SC}_1$ and this means exactly that $X$ is a (compact) ANR.
    We are left with showing $\dim(X)\leq 1$.

    Assume otherwise, i.e., assume $\dim(X)\geq 2$.
    Since $X$ is paracompact, we have $\locdim(X)=\dim(X)\geq 2$.
    This means there exists $x_0\in X$ such that $\dim(D)\geq 2$ for each closed neighborhood $D$ of $x_0$.
    For each $k$ consider $D_k:=\{y\in X\ :\ d(y,x_0)\leq 1/k\}$.
    This defines a decreasing sequence of closed neighborhoods around $x_0$ with $\dim(D_k)\geq 2$.

    It was noted in \cite[Proposition 3.1]{ChiDra2010} that a Peano space of dimension at least $2$ admits a topological
embedding\footnote{If $X,Y$ are spaces, then an injective map $i\colon X\to Y$ is called a topological embedding if the original topology of $X$ is the same as initial topology induced by the map $i$.
We usually consider a topologically embedded space as a subset with the subset topology.}
    of $S^1$.
    Indeed, a Peano space that contains no simple arc (i.e. in which $S^1$ cannot be embedded) is a dendrite, and therefore at most one-dimensional.
    It follows that there are embeddings $\varphi_k\colon S^1\hookrightarrow D_k\subset X$.
    Putting these together we get a map (not necessarily an embedding) $\varphi\colon Y\to X$ where $Y$ is the space of ''smaller and smaller circles'':
    \begin{align*}
        Y=\{(0,0)\}\cup\bigcup_{k\geq 1}S((1/2^k,0),1/(4\cdot 2^k))\subset\RR^2,
    \end{align*}
    where $S(x,r)$ is the circle of radius $r$ around the point $x$.
    We define $\varphi$ as $\varphi_k$ on the circle $S((1/k,0),1/3k)$.
    The map $\varphi\colon Y\to X$ induces a morphism $\varphi^\ast\colon C(X)\to C(Y)$.

    Next we construct a \Cs{} $B$ with a nested sequence of ideals $J_k\lhd B$, such that $C(Y)=B/\overline{\bigcup_k J_k}$ and $\varphi^\ast\colon C(X)\to C(Y)$ cannot be lifted to some $B/J_k$.
    Let $\mathcal{T}$ be the Toeplitz algebra and let $\mathcal{T}_1,\mathcal{T}_2,\ldots$ be a sequence of copies of the Toeplitz algebra, and set:
    \begin{align*}
        B &:=(\bigoplus_{k\in\NN} \mathcal{T}_k)^+ \\
        &=\{(b_1,b_2,\ldots)\in \prod_{k\geq 1}\mathcal{T} \text{ such that } (b_k)_k \text{ converges to a scalar multiple of } 1_\mathcal{T}\}.
    \end{align*}
    The algebras $\mathcal{T}_k$ come with ideals $\KK_k\lhd\mathcal{T}_k$ (each $\KK_k$ a copy of the algebra of compact operators $\KK$).
    Define ideals $J_k\lhd B$ as follows:
    \begin{align*}
        J_k &:=\KK_1\oplus\ldots\oplus\KK_k\oplus 0\oplus 0\oplus\ldots \\ %\quad \text{ ($k$ summands) } \\
        &=\{(b_1,\ldots,b_k,0,0,\ldots)\in B\ :\ b_i\in\KK_i\lhd\mathcal{T}_i\}.
    \end{align*}
    Note $B/J_k=C(S^1)\oplus\ldots_{(k)}\oplus C(S^1)\oplus (\bigoplus_{l\geq k+1}\mathcal{T}_l)^+$ ($k$ summands of $C(S^1)$).
    Also $J_k\subset J_{k+1}$ and $J:=\overline{\bigcup_k J_k}=\bigoplus_{k\in\NN}\KK_k$ and $B/J=(\bigoplus_{l\geq 1}C(S^1))^+\cong C(Y)$.

    The semiprojectivity of $C(X)$ gives a lift of $\varphi^\ast\colon C(X)\to C(Y)=B/J$ to some $B/J_k$.
    Consider the projection $\rho_{k+1}\colon B/J_k\to\mathcal{T}_{k+1}$ onto the (k+1)-th coordinate, and similarly $\varrho_{k+1}\colon B/J\to C(S^1)$.
    The composition $C(X)\to C(Y)\cong B/J\to C(S^1)$ is $\varphi_{k+1}^\ast$, the morphism induced by the inclusion $\varphi_{k+1}\colon S^1\hookrightarrow X$.
    Note that $\varphi_{k+1}^\ast$ is surjective since $\varphi_{k+1}$ is an inclusion.
    The situation is viewed in the following commutative diagram:
    \begin{center}
    	\makebox{
        \xymatrix{
        & & B/J_k \ar[d] \ar[r]^{\rho_{k+1}} & \mathcal{T}_{k+1} \ar[d] \\
        C(X) \ar@/_1pc/[rrr]_{\varphi_{k+1}^\ast} \ar[r]^>>{\varphi^\ast} \ar[urr] & C(Y) \ar[r]^{\cong} & B/J \ar[r]^{\varrho_{k+1}} & C(S^1)
        }}
    \end{center}
    The unitary $\id_{S^1}\in C(S^1)$ lifts under $\varphi_{k+1}^\ast$ to a normal element in $C(X)$, but it does not lift to a normal element in $\mathcal{T}_{k+1}$.
    This is a contradiction, and our assumption $\dim(X)\geq 2$ must be wrong.
    \qedhere
    \\
\end{proof}

%--------------------------------------------------------------------------------
\noindent
	It is well known that $C(D^2)$, the \Cs{} of continuous functions on the two-dimensional disc $D^2=\{(x,y)\in\RR^2\ :\ x^2+y^2\leq 1\}$, is not weakly semiprojective.
    For completeness we include the argument which is essentially taken from Loring \cite[17.1, p.131]{Lor1997}, see also \cite{Lor1995}.
    \\

%--------------------------------------------------------------------------------
\begin{prop}
\label{prop:S03:Disc_not_SP}
    $C(D^2)$ is not weakly semiprojective.
\end{prop}
\begin{proof}
    The $\ast$-homomorphisms from $C(D^2)$ to a \Cs{} $A$ are in natural one-one correspondence with normal contractions in $A$. Thus, statements about (weak) (semi-)projectivity of $C(D^2)$ correspond to statements about the (approximate) liftability of normal elements.
    For example, that $C(D^2)$ is projective would correspond to the (wrong) statement that normal elements lift from quotient \Cs s.
    To disprove weak semiprojectivity of $C(D^2)$ one uses a construction of operators that are approximately normal but do not lift in the required way due to an index obstruction.

    More precisely, define weighted shift operators $t_n$ on the separable Hilbert space $l^2$ (with basis $\xi_1,\xi_2,\ldots$) as follows:
    $$
        t_n(\xi_k)=\begin{cases}
            ((r+1)/2^{n-1})\xi_{k+1}   &\text{if } k=r2^{n+1}+s, 0\leq s<2^{n+1} \\
            \xi_{k+1}                 &\text{if } k\geq 4^n. \\
        \end{cases}
    $$
    Each $t_n$ is a finite-rank perturbation of the unilateral shift.
    Therefore the $t_n$ lie in the Toeplitz algebra $\Toep$ and have index $-1$.
    The construction of $t_n$ is made so that $\|t_n^\ast t_n-t_nt_n^\ast\|=1/2^{n-1}$.

    Consider the \Cs{} $B=\prod_\NN\Toep / \bigoplus_\NN\Toep$.
    The sequence $(t_1,t_2,\ldots)$ defines an element in $\prod_\NN\Toep$.
    Let $x=[(t_1,t_2,\ldots)]\in B$ be the equivalence class in $B$.
    Then $x$ is a normal element of $B$, and we let $\varphi\colon C(D^2)\to B$ be the corresponding morphism.
    We have the following lifting problem:

    \begin{center}\makebox{
    \xymatrix{
        & \prod_{k\geq N}\Toep_k \ar[d]^{\pi} \\
        C(D^2) \ar[r]^{\varphi} \ar@{..>}[ur]^{\bar{\varphi}} & {\prod_\NN\Toep / \bigoplus_\NN\Toep}  \\
    }}
    \end{center}

    Assume $C(D^2)$ is weakly semiprojective.
    Then the lifting problem can be solved, and $\bar{\varphi}$ defines a normal element $y=(y_N,y_{N+1},\ldots)$ in $\prod_{k\geq N}\Toep_k$.
    But the index of each $y_l$ is zero, while the index of each $t_l$ is $-1$, so that the norm-distance between $y_l$ and $t_l$ is at least one.
    Therefore the distance of $\pi(y)$ and $x$ is at least one, a contradiction.
    Thus, $C(D^2)$ is not weakly semiprojective.
    \qedhere
    \\
\end{proof}

%--------------------------------------------------------------------------------
\begin{rmk}[Spaces containing a two-dimensional disc]
\label{pargr:S03:spaces_containing_disc}
    We have seen above that $C(D^2)$ is not weakly semiprojective.
    Even more is true:
    Whenever a (compact, metric) space $X$ contains a two-dimensional disc, then $C(X)$ is not weakly semiprojective.
    This was noted by Loring, \cite{LorPrivat}.
    For completeness we include the argument:

    Let $D^2\subset X$ be a two-dimensional disc with inclusion map $i\colon D^2\to X$.
    Since $D^2$ is an absolute retract, there exists a retraction $r\colon X\to D$, i.e., $r\circ i=\id\colon D^2\to D^2$.
    Passing to \Cs s, we get induced momorphisms $i^\ast\colon C(X)\to C(D^2), r^\ast\colon C(D^2)\to C(X)$ such that $i^\ast\circ r^\ast$ is the identity on $C(D^2)$.
    Assume $C(X)$ is weakly semiprojective.
    Then any lifting problem for $C(D^2)$ could be solved as follows:
    Using the weak semiprojectivity of $C(X)$, the morphism $\varphi\circ i^\ast$ can be lifted.
    Then $\sigma\circ r^\ast$ is a lift for $\varphi=\varphi\circ i^\ast\circ r^\ast$.
    The situation is viewed in the following commutative diagram:
    \begin{center}\makebox{
    \xymatrix{
        & & & \prod_{k\geq N}B_k \ar[d]^{\pi} \\
        C(D^2) \ar[r]_{r^\ast}
            & C(X) \ar[r]_{i^\ast} \ar@{..>}[urr]^{\sigma}
            & C(D^2) \ar[r]_<<<<<{\varphi}  & {\prod_{k\geq 1} B_k / \bigoplus_{k\geq 1} B_k}  \\
    }}
    \end{center}
    This gives a contradiction, as we have shown above that $C(D^2)$ is not weakly semiprojective.

    However, that a space does not contain a two-dimensional disc is no guarantee that it has dimension at most one.
    These kind of questions are studied in continuum theory, and Bing, \cite{Bin1951}, gave examples of spaces of arbitrarily high dimension that are hereditarily
indecomposable\footnote{A continuum (i.e. compact, connected, metric space) is called decomposable if it can be written as the union of two proper subcontinua. Note that the union is not assumed to be disjoint. For example the interval $[0,1]$ is decomposable as it can be written as the union of $[0,1/2]$ and $[1/2,1]$. A continuum is called hereditarily indecomposable if none of its subcontinua is decomposable. See \cite{Nad1992} for further information.}, in particular they do not contain an arc or a copy of $D^2$.

    These pathologies cannot occur if we restrict to ''nicer'' spaces.
    For example, if a CW-complex does not contain a two-dimensional disc, then it has dimension at most one.
    What about ANRs?
    Bing and Borsuk, \cite{BinBor1964}, gave an example of a three-dimensional AR that does not contain a copy of $D^2$.
    The question for four-dimensional AR's is still open, i.e., it is unknown whether there exist high-dimensional AR's (or just ANRs) that do not contain a copy of $D^2$.

    The point we want to make clear is the following:
    To prove that an ANR is one-dimensional it is not enough to prove that it does not contain a copy of $D^2$.
    \\
\end{rmk}

%--------------------------------------------------------------------------------
\begin{rmk}[Spaces contained in ANRs of dimension $\geq 2$]
\label{pargr:S03:subspaces_ANR_dim2}
    Although an ANR $X$ of with $\dim(X)\geq 2$ might not contain a disc, one can show that it must contain (a copy of) one of the following three spaces:

    \begin{description}
        \item[Space 1]
        The space $Y_1$ of distinct ''smaller and smaller circles'' as considered in the proof of \ref{prop:S03:SP_gives_1D}, i.e.,
        $Y_1=\{(0,0)\}\cup\bigcup_{k\geq 1}S((1/2^k,0),1/(4\cdot 2^k))\subset\RR^2$.

        \item[Space 2]
        The Hawaiian earrings, i.e.,

        $Y_2=\bigcup_{k\geq 1}S((1/2^k,0),1/2^k)\subset\RR^2$.

        \item[Space 3]
        A variant of the Hawaiian earrings, where the circles do not just intersect in one point, but have a segment in common. It is homeomorphic to:

        $Y_3=\{(x,x),(x,-x)\ :\ x\in[0,1]\} \cup \bigcup_{k\geq 1} \{1/k\}\times[-1/k,1/k]\subset\RR^2$.
    \end{description}

    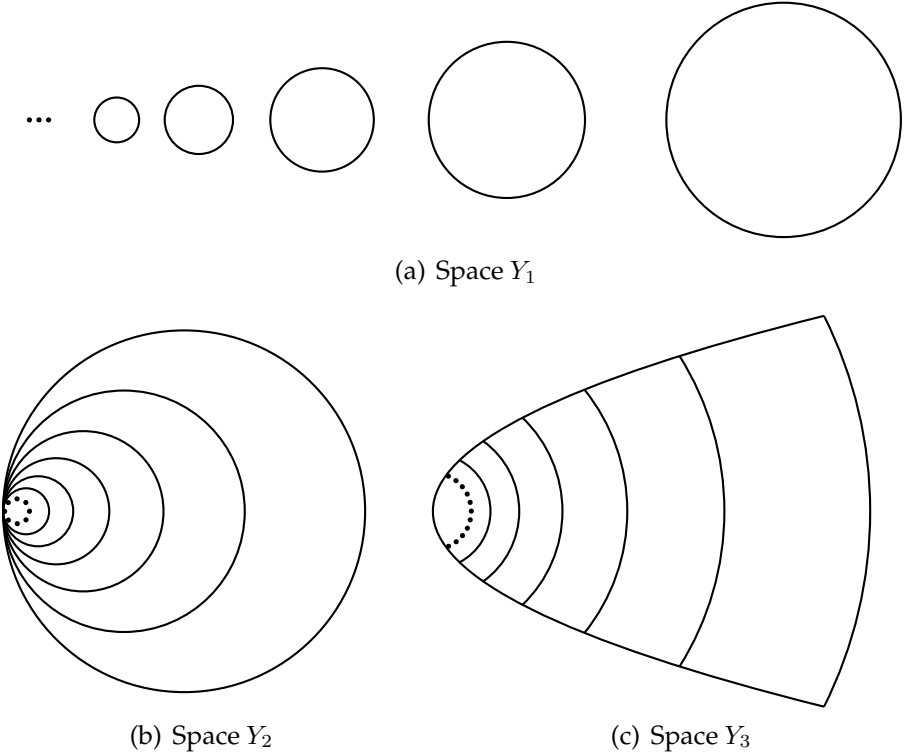
\begin{figure}[ht] \centering

    \subfigure[Space $Y_1$]{%\framebox{
        \psset{unit=0.013cm}
        \degrees
        \begin{pspicture}(-50,125)(900,-125)
            \psline[linestyle=dotted, linewidth=2pt, dotsep=2pt](-20,0)(5,0)
        	\pscircle(754, 0){121}
        	\pscircle(471, 0){81}
        	\pscircle(282, 0){54}

        	\pscircle(156, 0){36}
        	\pscircle(72, 0){24}
%        	\pscircle(16, 0){16}
        \end{pspicture}
    } %}

    \subfigure[Space $Y_2$]{%\framebox{
        \psset{unit=0.02cm}
        \begin{pspicture}(0,130)(260,-130)
        	\pscircle(121, 0){121}
        	\pscircle(81, 0){81}
        	\pscircle(54, 0){54}
        	\pscircle(36, 0){36}
        	\pscircle(24, 0){24}
        	\pscircle(16, 0){16}
        	\pscircle[linestyle=dotted, linewidth=2pt, dotsep=1.2pt](10, 0){10}
%          \psline[linestyle=dotted, linewidth=2.5pt, dotsep=1.2pt](0,0)(14,0)
        \end{pspicture}
    }%}
    \quad
    \subfigure[Space $Y_3$]{%\framebox{
        \psset{unit=1.3cm}
        \degrees
        \begin{pspicture}(0,2)(5,-2)
        	
        	\rput{-90}(2,2){
        	\parabola(0,2)(2,-2)
        	}

        	\psarc(0,0){4.472}{-26.565}{26.565}
        	\psarc(0,0){2.981}{-32.192}{32.192}
        	\psarc(0,0){1.988}{-38.776}{38.776}
        	\psarc(0,0){1.325}{-46.253}{46.253}
        	\psarc(0,0){0.883}{-54.334}{54.334}
        	\psarc(0,0){0.589}{-62.431}{62.431}
        	\psarc[linestyle=dotted, linewidth=2pt, dotsep=1.5pt](0,0){0.393}{-69.775}{69.775}
        \end{pspicture}
    }%}
    \caption{Spaces contained in high-dimensional ANRs}
    \end{figure}

    To prove this, one uses the same idea as in the proof of \ref{prop:S03:SP_gives_1D}:
    If $\dim(X)\geq 2$, then there exists a point $x_0$ where the local dimension is at least two.
    Then one can embed into $X$ a sequence of circles that get smaller and smaller and converge to $x_0$.
    Note that the circles may intersect or overlap.
    By passing to subspaces, we can get rid of ''unnecessary'' intersections and overlappings, and finally there are only three qualitatively different ways a bunch of ''smaller and smaller'' can look like.
    We skip the details.

    Note that none of the three spaces $Y_1,Y_2,Y_3$ are semiprojective.
    Further, no (compact, metric) space $X$ that contains a copy of $Y_1,Y_2$ or $Y_3$ can be semiprojective.
    One uses a similar argument as for an embedded $D^2$.
    Assume for some $k$ there is an inclusion $i\colon Y_k\hookrightarrow X$.
    Since $Y_k$ is not an AR, there will in general be no retraction onto it.

    Instead, choose an embedding $f\colon Y_k\hookrightarrow D^2$.
    This map can be extended a map $\tilde{f}\colon X\to D^2$ on all of $X$ since $D^2$ is an AR.
    \begin{center}\makebox{
    \xymatrix{
        D^2  \\
        Y_k \ar[r]_{i} \ar[u]^{f}
            & X \ar@{..>}[ul]_{\tilde{f}}
            & \\
    }}
    \end{center}

    If $C(X)$ is semiprojective, then any lifting problem as shown in the diagram below can be solved.
    However, using Toeplitz algebras as in \ref{prop:S03:SP_gives_1D} we see that the morphism $f^\ast=i^\ast\circ\tilde{f}^\ast\colon C(D^2)\to C(Y_k)$ is not semiprojective.
    \begin{center} \makebox{
    \xymatrix{
        & & & B/J_N \ar[d]^{\pi} \\
        C(D^2) \ar[r]_{\tilde{f}^\ast}
            & C(X) \ar[r]_{i^\ast} \ar@{..>}[urr]_{\sigma}
            & C(Y_k) \ar[r]_{\varphi}
            & B/\overline{\bigcup_{k\geq 1}J_k} \\
    }}
    \end{center}

    Finally let us note that the \Cs{s} $C(Y_1),C(Y_2)$ and $C(Y_3)$ are weakly semiprojective.
    \\
\end{rmk}

%################################################################################
%################################################################################
\section{Structure of compact, one-dimensional ANRs}
\label{sect:04:compact_ANR}

%--------------------------------------------------------------------------------
\noindent
	In this section we prove structural theorems about compact, one-dimensional
%\footnote{We say a space is one-dimensional if $\dim(X)\leq 1$. So, although it sounds weird, a one-dimensional space can also be zero-dimensional. It would probably be more precise to speak of ''at most one-dimensional'' space, however the usage of the term ''one-dimensional space'' is well established. }
    absolute neighborhood retracts (ANRs).
    The results are used in the next section to show that the \Cs{} of continuous functions on such a space is semiprojective.
    In section \ref{sect:06:Applications} we will study the structure on non-compact, one-dimensional ANRs.
	We start with some preparatory lemmas.
    By $\pi(X,x_0)$ we denote the fundamental group of $X$ based at $x_0\in X$.
    Statements about the fundamental group often do not depend on the basepoint, and then we will simply write $\pi(X)$ to mean that any (fixed) basepoint may be chosen.
	\\

%--------------------------------------------------------------------------------
\begin{lma}
\label{prop:S04::lma:homotope_path_to_piecewise_arc}
	Let $X$ be a Hausdorff space.
	Assume $X$ has a simply connected covering space.
	Then every path in $X$ is homotopic (relative endpoints) to a path that is piecewise arc.
\end{lma}
\begin{proof}
	Let $p\colon \widetilde{X}\to X$ be a simply connected, Hausdorff covering space.
	Let $\alpha\colon [0,1]\to X$ be a path, and let $\widetilde{\alpha}\colon [0,1]\to \widetilde{X}$ be a lift.
	Then the image of $\widetilde{\alpha}$ is a Peano continuum (i.e., a compact, connected, locally connected, metric space), and is therefore arcwise connected.
	Choose any arc $\beta\colon [0,1]\to \widetilde{X}$ from $\widetilde{\alpha}(0)$ to $\widetilde{\alpha}(1)$.
	The arc may of course be chosen within the image of $\widetilde{\alpha}$.
	Since $\widetilde{X}$ is simply connected, the paths $\widetilde{\alpha}$ and $\beta$ are homotopic (relative endpoints).
	Then $\alpha=p\circ\widetilde{\alpha}$ and $p\circ\beta$ are homotopic paths in $X$.
	
	Since $p$ is locally a homeomorphism, $p\circ\beta$ is piecewise arc, i.e., there exists a finite subdevision $0=t_0<t_1<\ldots<t_N=1$ such that each restriction $p\circ\beta_{|[t_j,t_{j+1}]}$ is an arc.
	\qedhere
	\\
\end{proof}

%--------------------------------------------------------------------------------
\begin{lma}
\label{prop:S04::lma:subgr_giving_fundGp}
	Let $X$ be a Hausdorff space, and $x_0\in X$.
	Assume $X$ has a simply connected covering space, and $\pi(X,x_0)$ is finitely generated.
	Then there exists a finite graph $Y\subset X$ with $x_0\in Y$ such that $\pi(Y,x_0)\to\pi(X,x_0)$ is surjective.
\end{lma}
\begin{proof}
	Choose a set of generators $g_1,\ldots,g_k$ for $\pi(X,x_0)$, represented by loops $\alpha_1,\ldots,\alpha_k\colon S^1\to X$.
	From the above lemma we can homotope each $\alpha_j$ to a loop $\beta_j$ that is piecewise arc.
	Then the image of each $\beta_j$ in $X$ is a finite graph.
	Consequently, also the union $Y:=\bigcup_j\im(\beta_j)$ is a finite graph (containing $x_0$).
	By construction each $g_j$ lies in the image of the natural map $\pi(Y,x_0)\to\pi(X,y_0)$.
	Therefore this map is surjective.	
	\qedhere
	\\
\end{proof}

%--------------------------------------------------------------------------------
\begin{rmk}
\label{pargr:S04:Existence_of_univ_cov_sp}
	Let $X$ be a connected, locally pathwise connected space.
	Then $X$ has a simply connected covering space (also called universal cover) if and only if $X$ is
    semilocally simply
connected\footnote{A space $X$ is called semilocally simply connected (sometimes also called locally relatively simply connected) if for each $x_0\in X$ there exists a neighborhood $U$ of $x_0$ such that $\pi(U,x_o)\to\pi(X,x_0)$ is zero.} (s.l.s.c.), see \cite[Theorem III.8.4, p.155]{Bre1993}.
	\\
\end{rmk}

%--------------------------------------------------------------------------------
\begin{prop}
\label{prop:S04:Subgr_fundGp_for_SLSC_Peano}
	Let $X$ be a s.l.s.c. Peano continuum and $x_0\in X$.
	Then there exists a finite graph $Y\subset X$ with $x_0\in Y$ such that $\pi(Y,x_0)\to\pi(X,x_0)$ is surjective.
\end{prop}
\begin{proof}
    Peano continua are connected and locally pathwise connected.
    Therefore, by the above remark \ref{pargr:S04:Existence_of_univ_cov_sp}, $X$ has a simply connected covering space.
	By \cite[Lemma 7.7]{CanCon2006}, $\pi(X,x_0)$ is finitely generated (even finitely presented).
	Now we may apply the above lemma \ref{prop:S04::lma:subgr_giving_fundGp}.
	\qedhere
	\\
\end{proof}

%--------------------------------------------------------------------------------
\begin{rmk}
\label{pargr:S04:subgraph_has_free_fundGp}
	The fundamental group of a finite graph is finitely generated (f.g.), free and abelian.
	Thus, the above map $\pi(Y,x_0)\to\pi(X,x_0)$ will in general not be injective.
	
	Even if $\pi(X,x_0)$ is f.g., free and abelian, the constructed map might not be injective.
	The reason is simply that the constructed graph could contain ''unnecessary'' loops
	(e.g. consider a circle embedded into a disc).
	However, by restricting to a subgraph one can get $\pi(Y,x_0)\to\pi(X,x_0)$ to be an isomorphism.
	
	Thus, if $X$ is a Hausdorff space that has a simply connected covering space, and $\pi(X,x_0)$ is finitely generated, free and abelian, then there exists a finite graph $Y\subset X$ such that $\pi(Y,x_0)\to\pi(X,x_0)$ is an isomorphism.	

	Let us consider a one-dimensional space $X$.
	This situation is special, since Cannon and Conner, \cite[Corollary 3.3]{CanCon2006}, have shown that an inclusion $Y\subset X$ of one-dimensional spaces induces an injective map on the fundamental group.
	Thus, we get the following:
	\\
\end{rmk}

%--------------------------------------------------------------------------------
\begin{prop}
\label{prop:S04:Subgr_fundGp_for_SLSC_1D}
	Let $X$ be a one-dimensional, Hausdorff space, and $x_0\in X$.
	Assume $X$ has a simply connected covering space, and $\pi(X,x_0)$ is finitely generated.
	Then there exists a finite graph $Y\subset X$ with $x_0\in Y$ such that $\pi(Y,x_0)\to\pi(X,x_0)$ is an isomorphism.
	\\
\end{prop}

%--------------------------------------------------------------------------------
\noindent
    Above we have studied, when there is a finite subgraph containing (up to homotopy) all loops of a space.
    We now turn to the question, when there is canonical such subgraph.
    It is clear that we can only hope for this to happen if the space is one-dimensional.

    We will use results from the master thesis of Meilstrup, \cite{Mei2005}, where also the following concept is introduced:
    A one-dimensional Peano continuum is called a \termDef{core continuum} if it contains no proper deformation retracts.
    \\

%--------------------------------------------------------------------------------
\begin{prop}[{see \cite[Corollary 2.6]{Mei2005}}]
\label{prop:S04:TFAE_core_continuum}

\noindent
	Let $X$ be a one-dimensional Peano continuum.
	Then the following are equivalent:
	\begin{enumerate}[label=(\arabic*)]
		\item
		$X$ is a core
		\item
		$X$ has no attached dendrites
        (an attached dendrite is a dendrite $C\subset X$ such that for some $y\in C$ there is a strong deformation retract $r\colon X\to(X\setminus C)\cup\{y\}$)
		\item
		every point of $X$ is on an essential loop that cannot be homotoped off it
        \item
        whenever $Y\subset X$ is a subset with $\pi(Y)\to\pi(X)$ surjective (hence bijective), then $Y=X$
	\end{enumerate}
\end{prop}
\begin{proof}
    The equivalence of (1),(2) and (3) is proved in \cite[Corollary 2.6]{Mei2005}.

    \impliesStep{3}{4}:
    Let $Y\subset X$ be a subset with $\pi(Y)\to\pi(X)$ surjective.
    Let $x\in X$ be any point.
    Then $x$ is on an essential loop, say $\alpha$, which cannot be homotoped off it.
    Since $[\alpha]\in\pi(Y,x)$ there is a loop $\beta$ with image in $Y$ that is homotopic to $\alpha$.
    Therefore $x\in Y$.

    \impliesStep{3}{4}:
    For any subset $Y$ that is a deformation retract of $X$ the map $\pi(Y)\to\pi(X)$ surjective.
    \qedhere
    \\
\end{proof}

%--------------------------------------------------------------------------------
\noindent
	To proceed further and prove that every one-dimensional Peano continuum contains a core we need the notion of reduced loop from \cite[Definition 3.8]{CanCon2006}.
	In fact, we will slighty generalize this to the notion of reduced path.
	This will help to simplify some proofs below.
	\\

%--------------------------------------------------------------------------------
\begin{defn}[{see \cite[Definition 3.8]{CanCon2006}}]
\label{defn:S04:reduced_path}

\noindent
		A non-constant path $\alpha\colon [0,1]\to X$ is called \termDef{reducible}, if there is an open arc $I=(s,t)\subset[0,1]$ such that $f(s)=f(t)$ and the loop $\alpha_{|[s,t]}$ based at $f(s)$ is nullhomotpic.
		A path is called \termDef{reduced} if it is not reducible.
		A constant path is also called reduced.
		\\
\end{defn}

%--------------------------------------------------------------------------------
\noindent
	By \cite[Theorem 3.9]{CanCon2006} every loop is homotopic to a reduced loop, and if the space is one-dimensional, then this reduced loop is even unique (up to reparametrization of $S^1$).
	The analogue for paths is proved in the same way.
	\\

%--------------------------------------------------------------------------------
\begin{prop}[{see \cite[Theorem 3.9]{CanCon2006}}]
\label{prop:S04:path_htpc_to_reduced_path}

\noindent
	Let $X$ be a space, and $\alpha\colon [0,1]\to X$ a path.
	Then $\alpha$ is homotopic (relative endpoints) to a reduced path $\beta\colon [0,1]\to X$ and we may assume the homotopy takes place inside the image of $\alpha$, so that also the image of $\beta$ lies inside the image of $\alpha$.
	If $X$ is one-dimensional, then the reduced path is unique up to reparametrizing of $[0,1]$.
	\\
\end{prop}

%--------------------------------------------------------------------------------
\begin{prop}[{see \cite[Theorem 2.4]{Mei2005}}]
\label{prop:S04:Existence_of_core}

\noindent
    Let $X$ be a non-contractible, one-dimensional Peano continuum.
    Then there exists a unique strong deformation retract $C\subset X$ that is a core continuum.
    We call it the core of $X$ and denote it by $\core(X)$.
    Further:
    \begin{enumerate}[label=(\arabic*)]
        \item
        $\core(X)$ is the smallest strong deformation retract of $X$
        \item
        $\core(X)$ is the smallest subset $Y\subset X$ such that the map $\pi(Y)\to\pi(X)$ is surjective
    \end{enumerate}
\end{prop}
\begin{proof}
    Let $\core(X)\subset X$ be the union of all essential, reduced loops in $X$.
    In the proof of \cite[Theorem 2.4]{Mei2005} it is shown that $\core(X)$ is a core continuum and a strong deformation retract of $X$.

    For every strong deformation retract $Y\subset X$ the map $\pi(Y)\to\pi(X)$ is surjective.
    Thus, to prove the two statements it is enough to show that $\core(X)$ is contained in every subset $Y\subset X$ such that the map $\pi(Y)\to\pi(X)$ is surjective.

    Let $Y\subset X$ be any subset such that the map $\pi(Y)\to\pi(X)$ is surjective, and let $\alpha$ be an essential, reduced loop in $X$.
    Then $\alpha$ is homotopic to a loop $\alpha'$ in $Y$.
    By the above remark the image of $\alpha'$ contains the image of $\alpha$.
    Thus, $Y$ contains all essential, reduced loops in $X$, and therefore $\core(X)\subset Y$.
    \qedhere
    \\
\end{proof}

%--------------------------------------------------------------------------------
\begin{rmk}
\label{pargr:S04:Core}
    If $X$ is a contractible, one-dimensional Peano continuum (i.e. a dendrite), then it can be contracted to any of its points.
    That is why $\core(X)$ is not defined in this situation.
    However, to simplify the following statements we will consider the core of a dendrite to be just any fixed point.

    If $X$ is a finite graph, then the core is obtained by successively removing all ''loose'' edges, i.e., vertices that are endpoints and the edge connecting the endpoint to the rest of the graph.
	\\
\end{rmk}

%--------------------------------------------------------------------------------
\noindent
    Next, we combine a bunch of known facts with some of our results to obtain a list of equivalent characterizations when a one-dimensional Peano continuum is an ANR.
    \\

%--------------------------------------------------------------------------------
\begin{thm}
\label{prop:S04:TFAE_1D_Peano_ANR}
	Let $X$ be a one-dimensional Peano continuum.
	Then the following are equivalent:
	\begin{enumerate}[label=(\arabic*)]
		\item
		$X$ is an absolute neighborhood retract (ANR)
		\item
		$X$ is locally contractible
		\item
		$X$ has a simply connected covering space
		\item
		$\pi(X)$ is finitely generated
		\item
		there exists a finite graph $Y\subset X$ such that $\pi(Y)\to\pi(X)$ is an isomorphism
		\item
		$\core(X)$ is a finite graph
	\end{enumerate}
\end{thm}
\begin{proof}
    \impliesStep{1}{2}:
    Every ANR is locally contractible, see \cite[V.2.3, p.101]{Bor1967}.

    \impliesStep{2}{3}:
    By the above remark \ref{pargr:S04:Existence_of_univ_cov_sp}.

    \impliesStep{3}{4}:
    By \cite[Lemma 7.7]{CanCon2006}.

    \impliesStep{4}{1}:
    This follows from \cite[V.13.6, p.138]{Bor1967}.

    ''(3)+(4) $\Rightarrow$ (5)'':
    Follows from \ref{prop:S04:Subgr_fundGp_for_SLSC_1D}.

    \impliesStep{5}{6}:
    By \ref{prop:S04:Existence_of_core} (2), $\core(X)\subset Y$.
    Then $\pi(\core(X))\to\pi(Y)$ is an isomorphism, and therefore $\core(X)=\core(Y)$.
    By the above remark \ref{pargr:S04:Core} the core of a finite graph is again a finite graph.

    \impliesStep{6}{4}:
    Follows since $\pi(\core(X))\to\pi(X)$ is bijective and the fundamental group of a finite graph is finitely generated.
    \qedhere
    \\
\end{proof}

%--------------------------------------------------------------------------------
\begin{rmk}
\label{pargr:S04:1D_Peano_AR_iff_core_pt}
	Let $X$ be a one-dimensional Peano continuum.
    In the same way as the above theorem \ref{prop:S04:TFAE_1D_Peano_ANR} one obtains that the following are equivalent:
	\begin{enumerate}[label=(\arabic*)]
		\item
		$X$ is an absolute retract (AR)
		\item
		$X$ is contractible
		\item
		$X$ is simply connected
		\item
		$\pi(X,x_0)$ is zero
        \item
		there exists a finite tree $Y\subset X$ such that $\pi(Y,x_0)\to\pi(X,x_0)$ is an isomorphism (for any $x_0\in Y$)
		\item
		$\core(X)$ is a point
	\end{enumerate}
    Note that $X$ is a dendrite if and only if it is a one-dimensional Peano continuum that satisfies one (or equivalently all) of the above conditions.
    \\
\end{rmk}

%--------------------------------------------------------------------------------
\noindent
	Let us proceed with the study of the internal structure of compact, one-dimensional ANRs.
	We will give a structure theorem which says that these spaces can be approximated by finite graphs in a nice way, namely from within.
	This generalzes a theorem from Nadler's book, \cite{Nad1992}, about the structure of dendrites (which are exactly the \emph{contractible} one-dimensional, compact ANRs).
	The point is that compact, one-dimensional ANRs can be approximated from within by finite graphs in exactly the same way as dendrites can be approximated by finite trees (which are exactly the contractible finite graphs).
	\\

%--------------------------------------------------------------------------------
\begin{lma}
\label{prop:S04::lma:First_pt_unique}
	Let $X$ be a one-dimensional Peano continuum, and $Y$ a subcontinuum with $\core(X)\subset Y$.
	For each $x\in X\setminus Y$ there is a unique point $r(x)\in Y$ such that $r(x)$ is a point of an arc in $X$ from $x$ to any point of $Y$.
\end{lma}
\begin{proof}
	This is the analogue of \cite[Lemma 10.24, p.175]{Nad1992}.
	We use ideas from the proof of \cite[Theorem 2.4]{Mei2005}.
	Let $X,Y$ be given, and $x\in X\setminus Y$.

	Pick some point $y\in Y$.
	Since $X$ is arc-connected, there exists an arc $\alpha\colon [0,1]\to X$ starting at $\alpha(0)=x$ and ending at $\alpha(1)=y$.
	Let $y_0=\alpha(\min\alpha^{-1}(Y))$, which is the first point in $Y$ of the arc (starting from $x$).
	Note that $y_0\in Y$ since $Y$ is closed.
	
	Assume there are two arcs $\alpha_1,\alpha_2\colon [0,1]\to X$ from $x$ to different points $y_1,y_2\in Y$ such that $\alpha_i([0,1))\subset X\setminus Y$.
	We show that this leads to a contradiction.
	Let $\beta$ be a reduced path in $Y$ from $y_1$ to $y_2$.
	Define
	\begin{align*}
		t_1	&:=\sup\{t\in[0,1]\ :\ \alpha_1(t)\in\im(\alpha_2)\} \\
		t_2	&:=\sup\{t\in[0,1]\ :\ \alpha_2(t)\in\im(\alpha_1)\},
	\end{align*}
	so that $x_0=\alpha_1(t_1)=\alpha_2(t_2)$ is the first point where the arcs $\alpha_1,\alpha_2$ meet (looking from $y_1$ and $y_2$).
	Connecting $(\alpha_1)_{|[t_1,1]}$ (from $x_0$ to $y_1$) with $\beta$ (from $y_1$ to $y_2$) and the inverse of $(\alpha_1)_{|[t_2,1]}$ (from $y_2$ to $x_0$), we get a reduced loop containing $x_0$ which contradicts $x_0\notin\core(X)\subset Y$.
    It follows that there exists a unique point $y\in Y$ with the desired properties.
	\qedhere
	\\
\end{proof}

%--------------------------------------------------------------------------------
\begin{defn}[{see \cite[Definition 10.26, p.176]{Nad1992}}]
\label{defn:S04:First_pt_map}

\noindent
	Let $X$ be a one-dimensional Peano continuum, and $Y$ a subcontinuum with $\core(X)\subset Y$.
	Define a map $r\colon X\to Y$ by letting $r(x)$ as in the lemma \ref{prop:S04::lma:First_pt_unique} above if $x\in X\setminus Y$, and $r(x)=x$ if $x\in Y$.
	This map is called the \termDef{first point map}.
	\\
\end{defn}

%--------------------------------------------------------------------------------
\noindent
	The first point map is continuous, and thus a retraction of $X$ onto $Y$.
	This is the analogue of \cite[Lemma 10.25, p.176]{Nad1992} and proved the same way.
	
	But more is true: As in the proof of \cite[Theorem 2.4]{Mei2005}, one can show that $Y$ is a strong deformation retract of $X$.
	\\

%--------------------------------------------------------------------------------
\begin{prop}
\label{prop:S04:First_pt_map_cts}
	Let $X$ be a one-dimensional Peano continuum, and $Y$ a subcontinuum with $\core(X)\subset Y$.
	Then the first point map is continuous.
	Further, there is a strong deformation retraction to the first point map.
\end{prop}
\begin{proof}
	Let $X,Y$ be given.
	As in the proof of \cite[Theorem 2.4]{Mei2005}, the complement $X\setminus Y$ consist of a collection of attached dendrites $\{C_i\}$.
	That means each $C_i\subset X$ is a dendrite such that $C_i\cap Y$ consists of exactly one point $y_i$ and such that there is a strong deformation retract $r_i\colon X\to(X\setminus C_i)\cup\{y_i\}$.
	Meilstrup shows that these strong deformation retracts can be assembled to give a strong deformation retract to the first point map $r$.
	\qedhere
	\\
\end{proof}

%--------------------------------------------------------------------------------
\begin{thm}
\label{prop:S04:Structure_1D_Peano}
	Let $X$ be a one-dimensional Peano continuum.
	Then there is a sequence $\{Y_k\}_{k=1}^\infty$ such that:
	\begin{enumerate}[label=(\arabic*)]
		\item
		each $Y_k$ is a subcontinuum of $X$
		\item
		$Y_k\subset Y_{k+1}$
		\item
		$\lim_k Y_k=X$
		\item
		$Y_1=\core(X)$ and for each $k$, $Y_{k+1}$ is obtained from $Y_k$ by attaching a line segment at a point, i.e., $\overline{Y_{k+1}\setminus Y_k}$ is an arc with an end point $p_k$ such that $\overline{Y_{k+1}\setminus Y_k}\cap Y_k=\{p_k\}$
		\item
		letting $r_k\colon X\to Y_k$ be the first point map for $Y_k$ we have that $\{r_k\}_{k=1}^\infty$ converges uniformly to the identity map on $X$
	\end{enumerate}
	If $X$ is also ANR, then all $Y_k$ are finite graphs.
	If $X$ is even contractible (i.e. is an AR, or equivalently a dendrite), then $\core(X)$ is just some point, and all $Y_k$ are finite trees.
\end{thm}
\begin{proof}
	This is the analogue of \cite[Lemma 10.24, p.175]{Nad1992}, and the proof goes through if we use our analoguous lemmas \ref{prop:S04::lma:First_pt_unique} and \ref{prop:S04:First_pt_map_cts}.
	\qedhere
	\\
\end{proof}

%################################################################################
%################################################################################
\section{The other implication of the main theorem: Sufficiency}
\label{sect:05:sufficiency}

%--------------------------------------------------------------------------------
\noindent
    For this implication we aim to mirror the approach of Chigogidze and Dranishnikov, \cite{ChiDra2010}.
    However we first show how to go from $C(X)$ being a universal \Cs{} to $C(Y)$ being one, where $Y$ is obtained from $X$ by attaching a line segment at one point.
    This step is not needed in \cite{ChiDra2010}, since they are able to give a general description of the generators and relations of the relevant spaces. We have not been able to find such generators and relations, and doing so might be of independent interest.
    \\

%--------------------------------------------------------------------------------
\begin{lma}
\label{ExtendUniversal}
    Suppose $X$ is a space, that $C(X) = C^* \langle \cG \mid \mathcal{R} \rangle$ and that $\{ \hat{g} \mid g \in \cG \}$ is a generating set of $C(X)$ that fulfills $\mathcal{R}$. Let $Y$ be the space formed from $X$ by attaching a line segment at a point $v$, and let $\lambda_g = \hat{g}(v)$. Then $C(Y) = C^* \langle \cG \cup \{ h \} \mid \mathcal{R}' \rangle$, where
    	\[
    		\mathcal{R}' = \mathcal{R} \cup \{ g h = \lambda_g h \text{ and } g h = h g \mid g \in \cG \} \cup \{ 0 \leq h \leq 1 \}.
    	\]
\end{lma}
\begin{proof}
    Extending the $\hat{g}$ to $Y$ by letting them be constant on the added line segment and letting $\hat{h}$ be the function that is zero on $X$ and grows linearly to one on the line segment (identifying it with $[0,1]$), shows that that there is a generating family in $C(Y)$ that fulfills $\mathcal{R}'$.

    We will use \cite[Lemma 3.2.2, p.26]{Lor1997} to show that $C(Y)$ is universal for $\mathcal{R}$.
    By this lemma, it suffices to show, that whenever we have a family $\{ T_g \mid g \in \cG \cup \{ h \} \}$ of operators, on some Hilbert space $H$, that fulfills $\mathcal{R}$ and $\{ T_g \mid g \in \cG \}' = \CC I$, then we can find a morphism from $C(Y)$ to $B(H)$ taking $\hat{g}$ to $T_g$ for all $g \in \cG \cup \{ h \}$.

    Suppose we have such operators.
    Since $C(X)$ is commutative and $\mathcal{R}'$ forces $h$ to commute with all the other generators, we have that $T_g = \mu_g I$ for some $\mu_g \in \CC$, for all $g \in \cG \cup \{ h \}$.
    We need to find a morphism from $C(Y)$ to $\CC$. There are two cases.
    \begin{itemize}[leftmargin=20pt, itemsep=5pt]
    	\item
        \textbf{Case 1:} $\mu_h = 0$:
        In this case we can find a morphism $\phi \colon C(X) \to \CC$ such that $\phi(\hat{g}) = \mu_g$ for all $g \in \cG$, since $C(X) = C^* \langle \cG \mid \mathcal{R} \rangle$.
        Then $\phi = \ev_u$ for some point $u \in X$.
        The morphism $\ev_u \colon C(Y) \to \CC$ maps $\hat{h} = 0$ and $\hat{g} = \mu_g$, and thus is the required morphism.
    	
    	\item
        \textbf{Case 2:} $\mu_h \neq 0$
        Since $0 \leq T_h \leq 1$, we have $0 < \mu_h \leq 1$.
        For $g \in \cG$ we have
    		\[
    			\mu_g \mu_h I = T_g T_h = \lambda_g \mu_h I.
    		\]
    		So since $\mu_h \neq 0$, we have $\mu_g = \lambda_g$ for all $g \in \cG$.
        Let us now identify the added line segment with $[0,1]$.
        The morphism $\ev_{\mu_h} \colon C(Y) \to \CC$, takes $\hat{h}$ to $\mu_h$ and $\hat{g}$ to $\lambda_g = \mu_g$.
        Hence it is the required morphism.
    \end{itemize}
    \qedhere
    \\
\end{proof}

%--------------------------------------------------------------------------------
\noindent
    We now provide a slightly altered (in both proof and statement) version of \cite[Proposition 4.1]{ChiDra2010}.
    \\

%--------------------------------------------------------------------------------
\begin{lma}
\label{AddLine}
    Suppose $X$ is a one-dimensional finite graph, that $C(X) = C^* \langle \cG \mid \mathcal{R} \rangle$, that $\{ \hat{g} \mid g \in \cG \}$ is a generating set of $C(X)$ that fulfills $\mathcal{R}$, and that $\cG$ is finite.
    Let $Y$ be the space formed from $X$ by attaching a line segment at a point $v$.
    Suppose we have a commutative square
    \begin{center}
    	\makebox{
    		\xymatrix{
    			C(X) \ar[d]_{\iota} \ar[r]^{\psi}	&	C \ar[d]^{\pi} \\
    			C(Y) \ar[r]_{\phi}					& C/J
    		}
    	 }
    \end{center}
    where $J$ is an ideal in the unital \Cs{} $C$, $\pi$ is the quotient morphism, $\psi$ and $\phi$ are unital morphisms, and $\iota$ is induced by the retraction from $Y$ onto $X$, i.e., $\iota$ takes a function in $C(X)$ to the function in $C(Y)$ given by
    \[
    	\iota(f)(x) = \left\{ \begin{array}{rl}
    													f(x), & x \in X, \\
    													f(v), & x \text{ is in the added line segment}
    												\end{array} \right..
    \]
    Then for every $\e > 0$ we can find a morphism $\chi \colon C(Y) \to C$ such that $\pi \circ \chi = \phi$ and $\| (\chi \circ \iota)(\hat{g}) - \psi(\hat{g}) \| \leq  \e$ for every $g \in \cG$.
\end{lma}
\begin{proof}
    Throughout the proof we use the notation of Lemma $\ref{ExtendUniversal}$.

    Let $\delta > 0$ be given.
    We will construct a $\delta$-representation $\{ d_g \mid g \in \cG \cup \{ h \} \}$ of $\mathcal{R}'$ in $C$ such that $\pi(d_g) = \phi(\iota(\hat{g}))$ for $g \in \cG$ and $\pi(d_h) = \phi(\hat{h})$.

    Let $q_\kappa \colon X \to X$ be the map that collapses the ball $B_{\kappa/2}(v)$, fixes $X \setminus B_{\kappa}(v)$, and extends linearly in between. Since there are only finitely many $\hat{g}$, we can find $\kappa_0$ such that $\| q_{\kappa_0}^*(\hat{g}) - \hat{g} \| \leq \delta / 2 $, where $q_\kappa^*$ is the morphism on $C(X)$ induced by $q_\kappa$. For simpler notation we let $q = q_{\kappa_0}$, and put $w_g = q^*(\hat{g})$ for all $g \in \cG$.

    Let $f_0$ be a positive function in $C(X)$ of norm $1$ that is zero on $X \setminus B_{\kappa_0 / 2}(v)$ and $1$ at $v$. Observe that if $f \in q^*(C_0(X \setminus \{v \}))$, then $f f_0 = 0$. Since $\hat{h} \leq \iota(f_0)$ and $\psi(f_0)$ is a lift of $\phi(\iota(f_0))$, we can, by \cite[Corollary 8.2.2, p.63]{Lor1997}, find a lift $\bar{h}$ of $\phi(\hat{h})$ such that $0 \leq \bar{h} \leq \psi(f_0)$. We now claim that $\{ \psi(\hat{g}) \mid g \in \cG\} \cup \{ \bar{h} \}$, is a $\delta$-representation of $\mathcal{R}$.

    Since the $\bar{g}$ fulfill the relations $\mathcal{R}$ and $\bar{h}$ is a positive contraction, we only need to check that $\psi(\hat{g})$ and $\bar{h}$ almost commute, and that $\psi(\hat{g}) \bar{h}$ is almost $\lambda_g \bar{h}$.

    First we note that since $0 \leq \bar{h} \leq \psi(f_0)$ for any $f \in q^*(C_0(X \setminus \{ v \}))$ we have
    \[
    	\| \psi(f) \bar{h}^{1/2} \|^2 = \| \psi(f) \bar{h} \psi(f)^* \| \leq \| \psi(f) \psi(f_0) \psi(f)^* \| = 0.
    \]
    Thus $\psi(f) \bar{h} = 0$. In particular we have
    \[
    	\psi(w_g - \lambda_g) \bar{h} = 0.
    \]
    Now we have
    \begin{align*}
    	\| \psi(\hat{g}) \bar{h} - \bar{h} \psi(\hat{g}) \| &= \| \psi(\hat{g}) \bar{h} - \psi(w_g - \lambda_g) \bar{h} - \bar{h} \psi(\hat{g}) + \bar{h} \psi(w_g - \lambda_g)  \| \\
    	&= \| \psi(\hat{g} - w_g) \bar{h} + \lambda_g \bar{h} - \bar{h}(\psi(\hat{g} - w_g)) - \lambda_g \bar{h} \| \\
    	&\leq \|\bar{h}\| (\| \psi(\hat{g} - w_g) \| + \|\psi(\hat{g} - w_g)\|) \\
    	&\leq 2 \| \hat{g} - w_g \| \leq 2 \cdot \delta /2 = \delta,
    \end{align*}
    for all $g \in \cG$. Likewise we have
    \begin{align*}
     \| \psi(\hat{g})\bar{h} - \lambda_g \bar{h} \| &= \| \psi(\hat{g}) \bar{h} - \lambda_g \bar{h} - \psi(w_g - \lambda_g) \bar{h} \| \\
      &= \| \psi(\hat{g} - w_g) \bar{h} + \lambda_g  \bar{h} - \lambda_g \bar{h} \| \\
      &= \| \psi(\hat{g} - w_g) \bar{h} \| \leq \| \hat{g} - w_g \| \leq \delta/2 \leq \delta,
    \end{align*}
    for all $g \in \cG$. So $\{ \psi(g) \mid g \in \cG \} \cup \{ \bar{h} \}$ is indeed a $\delta$-representation of $\mathcal{R}'$. Further we have that $\pi(\psi(\hat{g})) = \phi(\iota(\hat{g}))$ and that $\pi(\bar{h}) = \phi(h)$.

    Since $X$ is a one-dimensional finite graph, $Y$ is also a one-dimensional finite graph, so $C(Y)$ is semiprojective by \cite[Proposition 16.2.1, p.125]{Lor1997}.
    By \cite[Theorem 14.1.4, p.106]{Lor1997} the relations $\mathcal{R}'$ are then stable. So the fact that we can find a $\delta$-representation for all $\delta$ implies that we can find a morphism $\chi \colon C(Y) \to C$ such that $\pi \circ \chi = \phi$ and $\| \chi(\iota(\hat{g})) - \psi(\hat{g}) \| \leq \e$ for all $g \in \cG$.
    \qedhere
    \\
\end{proof}

%--------------------------------------------------------------------------------
\noindent
    We are now ready to show that some inductive limits have good lifting properties. In particular if we have an initial lift then we can lift all that follows.
    \\

%--------------------------------------------------------------------------------
\begin{prop}
\label{InductiveLimitProjective}
    Suppose that $X$ is a compact space such that $C(X)$ can be written as an inductive limit $\varinjlim_n C(Y_n) = C(X)$, where each $Y_n$ is a finite graph, $Y_{n+1}$ is just $Y_n$ with a line segments attached at a point (as in Lemma \ref{AddLine}), and the bonding morphisms $\iota_{n,n+1}\colon C(Y_n)\to C(Y_{n+1})$ are as the morphism in Lemma \ref{AddLine}, i.e., induced by retracting the attached interval to the attaching point.

    If there is a unital morphism $\phi \colon C(X) \to C/J$, where $J$ is an ideal in a unital \Cs{} $C$, and a unital morphism $\psi_1 \colon C(Y_1) \to C$ such that $\pi \circ \psi_1 = \phi \circ \iota_{1,\infty}$, then there is a unital morphism $\bar{\psi} \colon C(X) \to C$ such that $\pi \circ \bar{\psi} = \phi$.
\end{prop}
\begin{proof}
    We have the following situation:
    \begin{center}
    	\makebox{
        \xymatrix{
        & & C \ar[d]^{\pi} \\
        C(Y_1) \ar[r]_{\iota_{1,\infty}} \ar[urr]^{\psi_1}
            & C(X) \ar[r]_{\phi} \ar@{..>}[ur]_{\bar{\psi}}
            & C/J
        \\
        }}
    \end{center}
    As $Y_1$ is a finite graph, $C(Y_1)$ is finitely generated. Thus $C(Y_1)$ is a universal \Cs{} for some finite set of generators and relations, $C(Y_1) = C^* \langle \cG_1 \mid \mathcal{R}_1 \rangle $, say. In view of Lemma \ref{ExtendUniversal} we can now assume that $C(Y_n) = C^* \langle \cG_n \mid \mathcal{R}_n \rangle$, where $\cG_1 \subseteq \cG_2 \cdots$, and likewise for the $\mathcal{R}_n$. We also get from Lemma \ref{ExtendUniversal} that all the $\cG_n$ and $\mathcal{R}_n$ are finite.

    Since we are given $\psi_1$, we can, using Lemma \ref{AddLine} inductively, for any sequence of positive numbers $(\e_n)$ find morphisms $\psi_n \colon C(Y_n) \to C$ for each $n > 1$ such that $\pi \circ \psi_n = \phi \circ \iota_{n,\infty}$ and such that $\| \psi_n(\hat{g}) - \psi_{n-1}(\hat{g}) \| \leq \e_n$ for the generators $\hat{g}$ of $C(Y_n)$.

    We now wish to define new morphisms $\chi_n \colon C(Y_n) \to C$ such that $\pi \circ \chi_n = \phi \circ \iota_{n,\infty}$ and $\chi_{n+1}$ extends $\chi_n$. To this end we define, for each $n \in \NN$, elements $\{ \bar{g}_n \mid g \in \cG_n \}$, by
    \[
    	\bar{g}_n = \lim_k \psi_{n+k}(\hat{g}).
    \]
    The sequence $(\psi_{n+k}(\hat{g}))$ is Cauchy if $\sum \e_n < \infty$, so we will assume that.
    We claim that for any $n \in \NN$ the elements $\{ \bar{g}_n \mid g \in \cG_n \}$ in $C$ fulfill $\mathcal{R}_n$.
    By \cite[Lemma 13.2.3, p.103]{Lor1997} the set $\{ \bar{g}_n \mid g \in \cG_n \}$ is an $\e$-representation of $\mathcal{R}_n$ for all $\e > 0$ since $\{ \psi_{n+k}(\hat{g}) \mid g \in \cG_n \}$ is a representation of $\mathcal{R}_n$ for all $k$.
    Thus $\{ \bar{g}_n \mid g \in \cG_n \}$ is a representation of $\mathcal{R}_n$. Observe that if $m \geq n$, then $\bar{g}_m = \bar{g}_n$, since $\bar{g}_m$ is the limit of a tail of the sequence $\bar{g}_n$ is the limit of. Thus, we will drop the subscripts, and simply say that we have elements $\{ \bar{g} \mid g \in \cup \cG_n \}$ such that for any $n \in \NN$ the set $\{ \bar{g} \mid g \in \cG_n \}$ fulfills $\mathcal{R}_n$. Now we can define the $\chi_n$. We put $\chi_n(\hat{g}) = \bar{g}$, for $g \in \cG_n$, and this extends to a morphism since $C(Y_n) \cong C^* \langle \cG_n \mid \mathcal{R}_n \rangle$. We get $\chi_{n_1} \circ \iota_{n,n+1} = \chi_n$ and $\pi \circ \chi_n = \phi \circ \iota$ by universality, since it holds on generators.

    By the universal property of an inductive limit we get a morphism $\chi \colon C(X) \to C$ such that $\pi \circ \chi = \phi$.
    \qedhere
    \\
\end{proof}

%--------------------------------------------------------------------------------
\begin{rmk}
    Using the structure theorem for dendrites, \cite[Theorem 10.27, p.176]{Nad1992}, see \ref{prop:S04:Structure_1D_Peano}, and the above Proposition \ref{InductiveLimitProjective} we may deduce that for a dendrite $X$ the \Cs{} $C(X)$ is projective in $\mathcal{S}_1$ (the category of unital \Cs{s}, see \ref{pargr:S02:wSP}).
    Thus, we recover the implication \impliesStep{1}{2} of \cite[Theorem 4.3]{ChiDra2010}.

    To elaborate: Each dendrite $X$ can be approximated from within by finite trees, i.e., $C(X)\cong\varinjlim C(Y_k)$ where $Y_1$ is just a single point and the trees $Y_k$ are obtained by successive attaching of line segments.
    Since $C(Y_1)=\CC$ is projective in $\mathcal{S}_1$, we obtain from \ref{InductiveLimitProjective} that morphisms from $C(X)$ into a quotients can be lifted, i.e., $C(X)$ is projective in $\mathcal{S}_1$.
    \\
\end{rmk}

%--------------------------------------------------------------------------------
\noindent
    We are now ready to prove our main theorem:
    \\

%--------------------------------------------------------------------------------
\begin{proof}[Proof of theorem \ref{MainTheorem}]
    The implication ''$(I) \Rightarrow (II)$'' is Proposition \ref{prop:S03:SP_gives_1D}.

    Let us prove ''$(II) \Rightarrow (I)$'':
    So assume $X$ is a compact ANR with $\dim(X)\leq 1$.
    Note that $X$ can have at most finitely many components $X_i$.
    If we can show that each $C(X_i)$ is semiprojective, then $C(X)=\bigoplus_i C(X_i)$ will be semiprojective (since semiprojectivity is preserved by finite direct sums, see \cite[Theorem 14.2.1, p.110]{Lor1997}).
    So we may assume $X$ is connected.

    Then theorem \ref{prop:S04:Structure_1D_Peano} applies, and we may find an increasing sequence $Y_1\subset Y_2\subset\ldots\subset X$ of finite subgraphs such that:
	\begin{enumerate}[label=(\arabic*)]
		\item
		$\lim_k Y_k=X$, i.e., $\overline{\bigcup_kY_k}=X$
		\item
        $Y_{k+1}$ is obtained from $Y_k$ by attaching a line segment at a point
	\end{enumerate}
    Then $C(X)=\varinjlim_k C(Y_k)$ where each bonding morphism $\iota_{k,k+1}\colon C(Y_k)\to C(Y_{k+1})$ is induced by the retraction from $Y_{k+1}$ to $Y_k$ that contracts $Y_{k+1}\setminus Y_k$ to the point $\overline{Y_{k+1}\setminus Y_k}\cap Y_k$.
    Suppose now that we are given a unital \Cs{} $C$ and an increasing sequence of ideals $J_1\lhd J_2\lhd\ldots\lhd C$ and a unital morphism $\sigma\colon C(X)\to C/\overline{\bigcup_k J_k}$.
    We need to find a lift $\bar{\sigma}\colon C(X)\to C/J_l$ for some $l$.

    Consider the unital morphism $\sigma\circ\iota_{1,\infty}\colon C(Y_1)\to C/\overline{\bigcup_k J_k}$.
    By \cite[Proposition 16.2.1, p.125]{Lor1997}, the initial \Cs{} $C(Y_1)$ is semiprojective.
    Therefore, we can find an index $l$ and a unital morphism $\alpha\colon C(Y_1)\to C/J_l$ such that $\pi_l\circ\alpha=\sigma\circ\iota_{1,\infty}$.
    This is viewed in the following commutative diagram:
    \begin{center}
    	\makebox{
        \xymatrix{
        & & & & C \ar[d] \\
        & & & & C/J_l \ar[d]^{\pi_l} \\
        C(Y_1) \ar[r]_{\iota_{1,2}} \ar@{..>}[urrrr]^{\alpha}
            & C(Y_2) \ar[r]
            & \ldots \ar[r]
            & C(X) \ar[r]_{\sigma}  & C/\overline{\bigcup_k J_k}
        \\
        }}
    \end{center}

    Now we can apply \ref{InductiveLimitProjective} to find a unital morphism $\bar{\sigma}\colon C(X)\to C/J_l$ such that $\pi_l\circ\bar{\sigma}=\sigma$.
    This shows that $C(X)$ is semiprojective.
    \qedhere
    \\
\end{proof}

%################################################################################
%################################################################################
\section{Applications}
\label{sect:06:Applications}

\noindent
    In this section we give applications of our findings.
    First, we characterize semiprojectivity of non-unital, separable commutative \Cs{s}.
    Building on this, we are able to confirm a conjecture of Loring in the particular case of commutative \Cs{s}.
    Then, we will study the semiprojectivity of \Cs{s} of the form $C_0(X,M_k)$.
    Finally, we will give a partial solution to the problem when a commutative \Cs{} is weakly (semi-)projective.
    To keep this article short, we will omit most of the proofs in this sections.
    \\

    To characterize semiprojectivity of non-unital commutative \Cs{s} we have to study the structure of non-compact, one-dimensional ANRs.
    We are particularly interested in the one-point compactifications of such spaces.
    The motivation are the following results:
    If $X$ is a locally compact, Hausdorff space, then naturally $\widetilde{C_0(X)}\cong C(\cpctnPt{X})$, where $\cpctnPt{X}$ is the one-point comapctification of $X$.
    Further, a \Cs{} $A$ is semiprojective if and only if $\widetilde{A}$ is semiprojective.
    Thus, $C_0(X)$ is semiprojective if and only if $C(\cpctnPt{X})$ is semiprojective.
    By our main result \ref{MainTheorem} this happens precisely if $\cpctnPt{X}$ is a one-dimensional ANR.

    The following result gives a topological characterization of such spaces.
%    For a good introduction to the theory of ends we refer to the survey paper of Eilers, \cite{Eil1994}.
    We derive a characterization of semiprojectivity for non-unital, separable commutative \Cs{s}, see corollary \ref{prop:S06:SP_for_non-compact}.
    We also show that $\cpctnPt{X}$ is a one-dimensional ANR if and only if every finite-point
compactification\footnote{A compactification of a space $X$ is a pair $(Y,\iota_Y)$ where $Y$ is a compact space, $\iota\colon X\to Y$ is an embedding and $\iota(X)$ is dense in $Y$. Usually the embedding is understood and one denotes a compactification just by the space $Y$.
A compactification $\gamma(X)$ of $X$ is called a finite-point compactification if the remainder $\gamma(X)\setminus X$ is finite.}
    of $X$ is a one-dimensional ANR.
    Using this, we can confirm a conjecture about the semiprojective of extensions in the commutative case, see \ref{prop:S06:ideal_with_fd_quotient} and \ref{pargr:S06:Conj_SP_extension}.
	\\

%--------------------------------------------------------------------------------
\begin{thm}
\label{prop:S06:cpctfn_ANR}
	Let $X$ be a one-dimensional, locally compact, separable, metric ANR.
	Then the following are equivalent:
	\begin{enumerate}[label=(\arabic*)]
		\item
		The one-point compactification $\cpctnPt{X}$ is an ANR
		\item
		$X$ has only finitely many compact components and also only finitely many components $C\subset X$ such that $\cpctnPt{C}$ is not a dendrite
		\item
		Every finite-point compactification of $X$ is an ANR
		\item
		Some finite-point compactification of $X$ is an ANR
        \\
	\end{enumerate}
\end{thm}

%--------------------------------------------------------------------------------
\begin{cor}
\label{prop:S06:SP_for_non-compact}
    Let $X$ be a locally compact, separable, metric space.
    Then the following are equivalent:
    \begin{enumerate}[label=(\arabic*)]
	   \item $C_0(X)$ is semiprojective.
	   \item $X$ is a one-dimensional ANR that has only finitely many compact components, and $X$ has also only finitely many components $C\subset X$ such that $\cpctnPt{C}$ is not a dendrite
        \\
    \end{enumerate}
\end{cor}

%--------------------------------------------------------------------------------
\begin{cor}
\label{prop:S06:ideal_with_fd_quotient}
	Let $A$ be a separable, commutative \Cs{}, and $I\lhd A$ an ideal.
	Assume $A/I$ is finite-dimensional, i.e.. $A/I\cong\CC^k$ for some $k$.
	Then $A$ is semiprojective if and only if $I$ is semiprojective.
\end{cor}
\begin{proof}
    Let $A=C_0(X)$ for a locally compact, separable, metric space $X$.
    Then $I=C_0(Y)$ for an open subset $Y\subset X$.
    Since $A/I$ is finite-dimensional, $X\setminus Y$ is finite.
    It follows that also $\cpctnPt{X}\setminus Y$ is finite, and so the closure $\overline{Y}\subset\cpctnPt{X}$ is a finite-point compactification of $Y$.
    Set $F:=\cpctnPt{X}\setminus\overline{Y}$ (which is also finite).
    Note that $\overline{Y}\subset\cpctnPt{X}$ is a component, so that $\cpctnPt{X}=\overline{Y}\sqcup F$.
    It follows that $\cpctnPt{X}$ is an ANR if and only $\overline{Y}$ is.
    Then we argue as follows:

    \begin{tabular}[itemsep=5pt]{lll}
        & $A=C_0(X)$ is semiprojective \\

        $\Leftrightarrow$
            & $\widetilde{A}=C(\cpctnPt{X})$ is semiprojective \\

        $\Leftrightarrow$
            & $\cpctnPt{X}$ is a one-dimensional ANR
            & [ by theorem \ref{MainTheorem} ] \\

        $\Leftrightarrow$
            & $\overline{Y}\subset\cpctnPt{X}$ is a one-dimensional ANR
            & [ since $\cpctnPt{X}=\overline{Y}\sqcup F$] \\

    \multirow{2}{*}{$\Leftrightarrow$} &
    \multirow{2}{*}{$\cpctnPt{Y}$ is a one-dimensional ANR}
            & [by theorem \ref{prop:S06:cpctfn_ANR} since $\overline{Y}$ is a \\
            & & \ finite-point compactification of $Y$ ] \\

        $\Leftrightarrow$
            & $\widetilde{I}=C(\cpctnPt{Y})$ is semiprojective
            & [ by theorem \ref{MainTheorem} ] \\

        $\Leftrightarrow$
            & $I=C_0(Y)$ is semiprojective
    \end{tabular}

    \qedhere
    \\
\end{proof}

%--------------------------------------------------------------------------------
\begin{rmk}
\label{pargr:S06:Conj_SP_extension}
	Let $A$ be a separable \Cs{}, and $I\lhd A$ an ideal so that the quotient is finite-dimensional.
	We get a short exact sequence:
    \begin{center}
    	\makebox{
        \xymatrix{
        0\ar[r] & I \ar[r] & A \ar[r] & F \ar[r] & 0
        \\
        }}
    \end{center}
  It was conjectured by Loring and also by Blackadar, \cite[Conjecture 4.5]{Bla2004}, that in this situation $A$ is semiprojective if and only if $I$ is semiprojective.
  One implication was recently proven by Dominic Enders, \cite{EndPrivat}, who showed that semiprojectivity passes to ideals when the quotient is finite-dimensional.

  The converse implication is in general not even known for $F=\CC$.
  Our above result \ref{prop:S06:ideal_with_fd_quotient} confirms this conjecture in the case that $A$ is commutative.
  \\
\end{rmk}

%--------------------------------------------------------------------------------
\noindent
    Let us now study the semiprojectivity of \Cs{s} of the form $C_0(X,M_k)$.
    \\

%--------------------------------------------------------------------------------
\begin{lma} \label{eval}
    Let $X$ be a locally compact metric space and let $k\in\NN$.
    If $\phi \colon C_0(X,M_k) \to M_k$ is a morphism then there is a unitary $u \in M_k$ and a unique point $x \in \cpctnPt{X}$ such that
    \[
	\phi = Ad_u \circ \ev_x.
    \]
    \\
\end{lma}

%--------------------------------------------------------------------------------
\begin{prop}
\label{matrixAR}
    Let $X$ be a locally compact, separable, metric space and let $k \in \NN$.
    If $C_0(X,M_k)$ is projective, then $\cpctnPt{X}$ is an AR.
\end{prop}
\begin{proof}
    Suppose we are given a compact metric space $Y$ with an embedding $\iota\colon \cpctnPt{X}\to Y$. Dualizing and embedding $C_0(X)$ into $C(\cpctnPt{X})$, we get the following diagram
    \begin{center}
    	\makebox{
        \xymatrix{
            & C_0(Y) \ar[d]^{\iota_*} \\
    	   C_0(X) \ar[r] & C(\cpctnPt{X})
        }}
    \end{center}
    Tensoring everything by the $k$ by $k$ matrices $M_k$, we get
    \begin{center}
    	\makebox{
        \xymatrix{
    				& C_0(Y, M_k) \ar[d]^{(\iota_*)_k} \\
    	C_0(X,M_k) \ar[r]	& C(\cpctnPt{X}, M_k)
        }}
    \end{center}

    Since $C_0(X,M_k)$ is projective, there is a morphism $\psi \colon C_0(X,M_k) \to C_0(Y,M_k)$ such that $(\iota_*)_k \circ \psi$ is the inclusion of $C_0(X,M_k)$ into $C(\cpctnPt{X},M_k)$.

    For each $y \in Y$ lemma \ref{eval} tells us that the morphism $\ev_y  \circ \psi$, has the form $Ad_{u_y} \circ \ev_{x_y}$ for some unitary $u_y \in M_k$ and some unique $x_y \in \cpctnPt{X}$. Hence we can define a function $\lambda \colon Y \to \cpctnPt{X}$ such that
    \[
    	\ev_y \circ \psi = Ad_{u_y} \circ \ev_{\lambda(y)}.
    \]	
    This map $\lambda$ is continuous.

    For each $x \in\cpctnPt{X}$ we have the following commutative diagram
    \begin{center}
    	\makebox{
        \xymatrix{
    		& C_0(Y, M_k) \ar[d]^{(\iota_*)_k} \ar[r]^>>>>{\ev_{\iota(x)}}	& M_k \ar@{=}[d] \\
    	C_0(X,M_k) \ar[r] \ar[ur]^{\psi}	& C(\cpctnPt{X}, M_k) \ar[r]^>>>>{\ev_x}				 & M_k					
        }}
    \end{center}
    From this diagram, it follows that if $x \in \cpctnPt{X}$ then
    \[
    	Ad_{u_{\iota(x)}} \circ \ev_{\lambda(\iota(x))} = \ev_{\iota(x)} \circ \psi = \ev_x \circ (\iota_*)_k \circ \psi = \ev_x.
    \]	
    So for any function $g \in C_0(X, M_k)$ we get
    \[
    	\ev_{\lambda(\iota(x))} \begin{pmatrix} g &  & \\ & \ddots & \\ & & g \end{pmatrix} = (Ad_{u_{\iota(x)}} \circ \ev_{\lambda(\iota(x))}) \begin{pmatrix} g &  & \\ & \ddots & \\ & & g \end{pmatrix} = \ev_x \begin{pmatrix} g &  & \\ & \ddots & \\ & & g \end{pmatrix}.
    \]
    Hence we must have $\lambda(\iota(x)) = x$.

    All in all, we have found a continuous map $\lambda \colon Y \to \cpctnPt{X}$ such that $\lambda \circ \iota = \id$, i.e., the embedded space $\cpctnPt{X}\subset Y$ is a retract.
    As the embedding was arbitrary, $\cpctnPt{X}$ is an AR.
    \qedhere
    \\
\end{proof}

%--------------------------------------------------------------------------------
\noindent
    The proof can be modified to show:
    \\

%--------------------------------------------------------------------------------
\begin{prop}
\label{matrixANR}
    Let $X$ be a locally compact, separable, metric space and let $k \in \NN$.
    If $C_0(X,M_k)$ is semiprojective, then $\cpctnPt{X}$ is an ANR.
    \\
\end{prop}

%--------------------------------------------------------------------------------
\noindent
    Using the idea of the proof of \ref{prop:S03:SP_gives_1D} one can show the following:
    \\

%--------------------------------------------------------------------------------
\begin{prop}
\label{matrixDimension}
    Let $X$ be a locally compact, separable, metric space, and $k\in\NN$.
    If $C_0(X, M_k)$ is semiprojective, then $\dim(X) \leq 1$.
    \\
\end{prop}

%--------------------------------------------------------------------------------
\begin{cor}
\label{matrixMain}
    Let $A$ be a separable, commutative \Cs{}, and $k\in\NN$.
    If $A\otimes M_k$ is projective, then so is $A$.
    Analogously, if $A\otimes M_k$ is semiprojective, then so is $A$.
\end{cor}
\begin{proof}
    Let $A=C_0(X)$ for a locally compact, separable, metric space $X$.

    First, assume $A\otimes M_k$ is semiprojective.
    By proposition \ref{matrixDimension}, $\dim(X)\leq 1$.
    This implies that the dimension of $\cpctnPt{X}$ is at most one.
    By proposition \ref{matrixANR}, $\cpctnPt{X}$ is an ANR.
    Then our main theorem \ref{MainTheorem} shows that $C(\cpctnPt{X})$ is semiprojective.
    Since $C(\cpctnPt{X})$ is the unitization of $C_0(X)$, we also have that $C_0(X)$ is semiprojective.

    Assume now that $A\otimes M_k$ is projective.
    It follows that $A$ cannot be unital, for otherwise $A\otimes M_k$ would be unital and that is impossible for projective \Cs{s}.
    As in the semiprojective case we deduce $\dim(\cpctnPt{X})\leq 1$.
    By \ref{matrixAR}, $\cpctnPt{X}$ is an AR.
    It follows from \cite[Theorem 4.3]{ChiDra2010}, see also \ref{summaryThm}, the $C(\cpctnPt{X})$ is projective in $\mathcal{S}_1$.
    It follows that $C_0(X)$ is projective, see \ref{pargr:S02:wSP}.
    \qedhere
    \\
\end{proof}

%--------------------------------------------------------------------------------
\noindent
    We will now turn to the question, when a unital, commutative \Cs{} is weakly (semi-)projective in $\mathcal{S}_1$.
    The analogue of a weakly (semi-)projective \Cs{} in the commutative world is an approximative absolute (neighborhood) retract (abbreviated by AAR and AANR).
    As mentioned in \ref{pargr:S02:Connection}, if $C(X)$ is weakly (semi-)projective, then $X$ is AA(N)R.
    We will show below, that for one-dimensional spaces the converse is also true.
    \\

%--------------------------------------------------------------------------------
\begin{pargr}
\label{pargr:S07:conditions_AANR}
    Let $X$ be a compact, metric space.
    Consider the following conditions:
    \begin{enumerate}[label=(\arabic*)]
        \item
        for each $\varepsilon>0$ there exists a map $f\colon X\to Y\subset X$ such that $Y$ is an AR (an ANR), and $d(f)\leq\varepsilon$
        \item
        $X$ is an AAR (an AANR)
    \end{enumerate}

    Here, by $d(f)<\e$ we mean that the distance of $x$ and $f(x)$ is less than $\e$ for all $x\in X$, i.e., $d(x,f(x))<\e$ for all $x\in X$.
    The first condition means that $X$ can be approximated from within by ARs (by ANRs).
    As shown by Clapp, \cite[Theorem 2.3]{Cla1971}, see also \cite[Proposition 2.2(a)]{ChaPra2005}, the implication \impliesStep{1}{2} holds in general.

    It was asked by Charatonik and Prajs, \cite[Question 5.3]{ChaPra2005}, whether the converse also holds (at least for continua).
    They showed that this is indeed the case for hereditarily unicoherent continua, \cite[Observation 5.4]{ChaPra2005}.
    In theorem \ref{prop:S07:TFAE_1D_AANR} below we show that the two conditions are also equivalent for one-dimensional, compact, metric spaces.
    \\
\end{pargr}

\noindent
    The following is a standard result from continuum theory:
    \\

%--------------------------------------------------------------------------------
\begin{prop}
\label{prop:S07:1D_Peano_inner_approx_by_graphs}
    Let $X$ be a one-dimensional Peano continuum, and $\e>0$.
    Then there exists a finite subgraph $Y\subset X$ and a surjective map $f\colon X\to Y\subset X$ such that $d(f)<\e$.
    \\
\end{prop}

%--------------------------------------------------------------------------------
\begin{cor}
    Every one-dimensional Peano continuum is an AANR.
\end{cor}
\begin{proof}
    Let $X$ be a one-dimensional Peano continuum.
    By \ref{prop:S07:1D_Peano_inner_approx_by_graphs}, $X$ can be approximated from within by finite subgraphs.
    A finite graph is an ANR.
    It follows from \cite[Theorem 2.3]{Cla1971}, see \ref{pargr:S07:conditions_AANR}, that $X$ is an AANR.
    \qedhere
    \\
\end{proof}

\noindent
    The following Lemma is a direct translation of \cite[Lemma 5.5]{Lor2009} to the commutative setting.
    \\

%--------------------------------------------------------------------------------
\begin{lma}[{see \cite[Lemma 5.5]{Lor2009}}]

\noindent
    Let $X$ be an compact AAR, and $D$ any ANR.
    Then every map $f\colon X\to D$ is inessential, i.e., homotopic to a constant map.
    \\
\end{lma}

%--------------------------------------------------------------------------------
\begin{cor}
\label{prop:S07:1D_AAR_tree-like}
    Every one-dimensional, compact AAR is tree-like.
\end{cor}
\begin{proof}
    Let $X$ be a one-dimensional, compact AAR.
    Then $X$ is connected and thus a continuum.
    In \cite[Theorem 1]{CasCha1960} tree-like continua are characterized as one-dimensional continua such that every map into a finite graph is inessential.
    Thus, we need to show that every map from $X$ into a finite graph is inessential.
    This follows from the above Lemma since every finite graph is an ANR.
    \qedhere
    \\
\end{proof}

%--------------------------------------------------------------------------------
\begin{thm}
\label{prop:S07:TFAE_1D_AANR}
    Let $X$ be a one-dimensional, compact, metric space.
    Then the following are equivalent:
    \begin{enumerate}[label=(\arabic*)]
        \item
        for each $\varepsilon>0$ there exists a map $f\colon X\to Y\subset X$ such that $Y$ is a finite tree (a finite graph),
        and $d(f)\leq\varepsilon$
        \item
        for each $\varepsilon>0$ there exists a map $f\colon X\to Y\subset X$ such that $Y$ is an AR (an ANR),
        and $d(f)\leq\varepsilon$
        \item
        $X$ is an AAR (an AANR)
    \end{enumerate}
    Moreover, in $(1)$ and $(2)$ the map $f$ may be assumed to be surjective.
\end{thm}
\begin{proof}
    \impliesStep{1}{2} is clear, and \impliesStep{2}{3} follows from \cite[Theorem 2.3]{Cla1971}, see \ref{pargr:S07:conditions_AANR}.

    \impliesStep{3}{1}:
    It was shown by Clapp, \cite[Theorem 4.5]{Cla1971}, that for each embedding of a compact AANR $X$ in the Hilbert cube $Q$ and $\delta>0$ there exists a compact polyhedron $P\subset Q$ with maps $f\colon X\to P$ and $g\colon P\to X$ such that $d(f)<\delta$ and $d(g)<\delta$.
    Note that $g$ maps each component of $P$ onto a Peano subcontinuum of $X$.
    Thus, the image $Y:=g(P)\subset X$ is a finite union of Peano subcontinua.
    Moreover, the map $g\circ f:X\to Y\subset X$ satisfies $d(f)<2\delta$.

    Assume $X$ is a one-dimensional, compact AANR and fix some $\e>0$.
    We apply the result of Clapp for $\delta=\e/4$ and obtain a compact subspace $Y\subset X$ that is the (disjoint) union of finitely many Peano continua, together with a surjective map $f\colon X\to Y$ such that $d(f)<\e/2$.
    Since $Y\subset X$ is closed, $\dim(Y)\leq\dim(X)\leq 1$.
    Applying \ref{prop:S07:1D_Peano_inner_approx_by_graphs} to each component of $Y$ and $\e/2$ we obtain a finite subgraph $Z\subset Y$ and a surjective map $g\colon Y\to Z$ such that $d(g)<\e/2$.

    We may consider $Z$ as a finite subgraph of $X$.
    The map $h:=g\circ f\colon X\to Z\subset X$ is surjective and satisfies $d(h)<\e$.
    So we have shown the implication for the case that $X$ is AANR.

    Assume additionally that $X$ is an AAR.
    We have already shown that $X$ can be approximated from within by finite subgraphs.
    We need to show that the same is true with finite trees.

    By \ref{prop:S07:1D_AAR_tree-like}, $X$ is tree-like.
    By \cite[2.2 and 2.3]{Lel1976}, every tree-like continuum is hereditarily unicoherent.
%\footnote{A continuum $X$ is called unicoherent if for each two subcontinua $Y_1,Y_2\subset X$ with $X=Y_1\cup Y_2$ the intersection $Y_1\cap Y_2$ is a continuum (i.e. connected).
%    A continuum is called hereditarily unicoherent if all its subcontinua are unicoherent.}.
    A coherent finite graph is a finite tree.
    It follows that \emph{every} finite subgraph $Z\subset X$ is a finite tree, and so $X$ can be approximated from within by finite subgraphs which automatically are finite trees.
    \qedhere
    \\
\end{proof}

%--------------------------------------------------------------------------------
\begin{cor}
\label{prop:S07:1D_AANR_implies_wSP}
    Let $X$ be a compact, metric space.
    Then the following implications hold:
    \begin{enumerate}[label=(\arabic*)]
	   \item
        If $X$ is an AANR and $\dim(X)\leq 1$, then $C(X)$ is weakly semiprojective $\mathcal{S}_1$.
	   \item
        If $X$ is an AAR and $\dim(X)\leq 1$, then $C(X)$ is weakly projective in $\mathcal{S}_1$.
    \end{enumerate}
\end{cor}
\begin{proof}
    Let $X$ be a one-dimensional, compact AAR (AANR).
    By \ref{prop:S07:TFAE_1D_AANR}, $X$ can be approximated from within by finite trees (finite graphs), i.e., for each $n\geq 1$ there exists a finite tree (graph) $Y_n\subset X$ and a surjective map $f_n\colon X\to Y_n$ with $d(f_n)<1/n$.
    We want to use \cite[Theorem 4.7]{Lor2009} to show $C(X)$ is weakly (semi-)projective in $\mathcal{S}_1$.

    The surjective maps $f_n$ induce injective morphisms $f_n^\ast\colon C(Y_n)\to C(X)$.
    Consider also the inclusion map $\iota_n\colon Y_n\hookrightarrow X$ and the dual morphism $\iota_n^\ast\colon C(X)\to C(Y_n)$.
    Set $\theta_n:=f_n^\ast\circ\iota_n^\ast\colon C(X)\to C(X)$.

    Since $d(f_n)$ tends to zero, the morphisms $\theta_n$ converge (pointwise) to the identity morphism.
    Further, the image of $\theta_n$ is equal to the image of $f_n^\ast$, and therefore isomorphic to $C(Y_n)$.

    As shown by Loring, \cite[Proposition 16.2.1, p.125]{Lor1997}, $C(Y)$ is semiprojective (in $\mathcal{S}_1$) if $Y$ is a finite graph.
    Similarly, $C(Y)$ is projective in $\mathcal{S}_1$ if $Y$ is a finite tree $Y$ (see also \cite{ChiDra2010}).
    Now, it follows from \cite[Theorem 4.7]{Lor2009} (and the analogous result for weakly semiprojective \Cs{s}) that $C(X)$ is weakly (semi-)projective in $\mathcal{S}_1$.
    \qedhere
    \\
\end{proof}

%--------------------------------------------------------------------------------
\begin{rmk}
    We remark that the converse implications of \ref{prop:S07:1D_AANR_implies_wSP} also hold.
    As explained in \ref{pargr:S02:Connection}, if $C(X)$ is weakly (semi-)projective in $\mathcal{S}_1$, then $X$ is necessarily an approximative absolute (neighborhood) retract.
    The dimension condition was recently shown by Enders, \cite{EndPrivat}.

    Thus, $C(X)$ is (weakly) (semi-)projective in $\mathcal{S}_1$ if and only if $X$ is a compact (approximative) absolute (neighborhood) retract with $\dim(X)\leq 1$.
    \\
\end{rmk}

\section*{Acknowledgments}

\noindent
    We thank Dominic Enders for his comments and inspiring suggestions that helped to improve some of the results in section \ref{sect:06:Applications}.

    We thank S{\o}ren Eilers for his valuable comments on the first draft of this paper.

\end{document}